\documentclass[final,onefignum,onetabnum]{siamart190516}

\pdfoutput=1 

\usepackage{lipsum}
\usepackage{amsfonts}
\usepackage{graphicx}
\usepackage{epstopdf}
\usepackage{algorithmic}
\usepackage{amsopn}

\usepackage{amssymb} 
\usepackage{subcaption} 
\usepackage{booktabs} 
\usepackage[italicdiff]{physics} 
\usepackage{mathtools} 
\usepackage{bbm} 
\let\originalleft\left 
\let\originalright\right
\renewcommand{\left}{\mathopen{}\mathclose\bgroup\originalleft}
\renewcommand{\right}{\aftergroup\egroup\originalright}
\newcommand{\one}{\mathbbm{1}} 
\newcommand{\E}{\mathbb{E}}
\newcommand{\N}{\mathbb{N}}
\newcommand{\Z}{\mathbb{Z}}
\newcommand{\R}{\mathbb{R}}
\newcommand{\cP}{\mathcal{P}}
\newcommand{\cA}{\mathcal{A}}
\newcommand{\cL}{\mathcal{L}}
\newcommand{\cM}{\mathcal{M}}
\newcommand{\cF}{\mathcal{F}}
\newcommand{\cH}{\mathcal{H}}
\newcommand{\cX}{\mathcal{X}}
\newcommand{\cY}{\mathcal{Y}}
\newcommand{\al}{\alpha}
\newcommand{\ep}{\varepsilon}
\DeclarePairedDelimiterX{\iptemp}[2]{\langle}{\rangle}{#1, #2}
\renewcommand{\ip}{\iptemp}
\DeclarePairedDelimiterX{\normtemp}[1]{\lVert}{\rVert}{#1}
\renewcommand{\norm}{\normtemp}
\DeclarePairedDelimiterX{\abstemp}[1]{\lvert}{\rvert}{#1}
\renewcommand{\abs}{\abstemp}
\newcommand{\ON}[1]{\operatorname{#1}}
\newcommand{\lap}{\Delta}
\DeclareMathOperator{\id}{Id} 
\DeclareMathOperator*{\argmin}{argmin}
\DeclareMathOperator{\im}{im} 
\renewcommand{\ker}{\ON{ker}} 
\newcommand{\spanvec}{\ON{span}} 
\newcommand{\xqed}[1]{%
	\leavevmode\unskip\penalty9999 \hbox{}\nobreak\hfill
	\quad\hbox{\ensuremath{#1}}}
\newcommand{\Endrmk}{\xqed{\lozenge}} 
\newcommand{\Fd}{F^{\dagger}} 
\newcommand{\defby}{\coloneqq} 
\newcommand{\diid}{\stackrel{\mathrm{iid}}{\sim}} 
\def\qfa{\quad\mbox{for all}\quad} 

\ifpdf
\DeclareGraphicsExtensions{.eps,.pdf,.png,.jpg}
\else
\DeclareGraphicsExtensions{.eps}
\fi
\graphicspath{ {figures/} }

\newsiamremark{remark}{Remark}
\newsiamremark{hypothesis}{Hypothesis}
\crefname{hypothesis}{Hypothesis}{Hypotheses}
\newsiamthm{claim}{Claim}
\newsiamremark{example}{Example}
\newsiamremark{notation}{Notation}
\newsiamthm{result}{Result}

\headers{The Random Feature Model on Banach Space}{N. H. Nelsen and A. M. Stuart}

\title{The Random Feature Model for Input-Output Maps between Banach Spaces\thanks{
		Submitted to the editors May 20, 2020; accepted for publication (in revised form) May 20, 2021.
		\funding{NHN is supported by the National Science Foundation (NSF) Graduate Research Fellowship Program under award DGE-1745301. AMS is supported by NSF (award DMS-1818977) and by the Office of Naval Research (ONR) (award N00014-17-1-2079). Both authors are supported by NSF (award AGS-1835860) and ONR (award N00014-19-1-2408).}
	}
}
\author{Nicholas H. Nelsen\thanks{Division of Engineering and Applied Science, California Institute of Technology, Pasadena, CA 91125, USA (\email{nnelsen@caltech.edu}, \email{astuart@caltech.edu}).}
	\and Andrew M. Stuart\footnotemark[2]
}

\ifpdf
\hypersetup{
	pdftitle={The Random Feature Model on Banach Space},
	pdfauthor={N. H. Nelsen and A. M. Stuart}
}
\fi

\begin{document}
	
\maketitle

\begin{abstract}
	Well known to the machine learning community, the random feature model is a parametric approximation to kernel interpolation or regression methods. It is typically used to approximate functions mapping a finite-dimensional input space to the real line. In this paper, we instead propose a methodology for use of the random feature model as a data-driven surrogate for operators that map an input Banach space to an output Banach space. Although the methodology is quite general, we consider operators defined by partial differential equations (PDEs); here, the inputs and outputs are themselves functions, with the input parameters being functions required to specify the problem, such as initial data or coefficients, and the outputs being solutions of the problem. Upon discretization, the model inherits several desirable attributes from this infinite-dimensional viewpoint, including mesh-invariant approximation error with respect to the true PDE solution map and the capability to be trained at one mesh resolution and then deployed at different mesh resolutions. We view the random feature model as a non-intrusive data-driven emulator, provide a mathematical framework for its interpretation, and demonstrate its ability to efficiently and accurately approximate the nonlinear parameter-to-solution maps of two prototypical PDEs arising in physical science and engineering applications: viscous Burgers' equation and a variable coefficient elliptic equation.
\end{abstract}



\section{Introduction}
\label{sec:intro}
The \emph{random feature model}, an architecture for the data-driven approximation of maps between finite-dimensional spaces, was formalized in 
\cite{rahimi2008random, rahimi2008uniform, rahimi2008weighted}, building on earlier precursors in~\cite{barron1993universal,neal1996priors,williams1997computing}. The goal of this paper is to extend the random feature model
to a methodology for the data-driven approximation of maps between infinite-dimensional spaces. Canonical examples of such maps include the semigroup generated by a time-dependent partial differential equation (PDE) mapping the initial condition (an input parameter) to the solution at a later time and the operator mapping a coefficient function (an input parameter) appearing in a PDE to its solution. Obtaining efficient and potentially low-dimensional representations of PDE solution maps is not only conceptually interesting, but also practically useful. Many applications in science and engineering require repeated evaluations of a complex and expensive forward model for different configurations of a system parameter. The model often represents a discretized PDE and the parameter, serving as input to the model, often represents a high-dimensional discretized quantity such as an initial condition or uncertain coefficient field. These \emph{outer loop} applications commonly arise in inverse problems or uncertainty quantification tasks that involve control, optimization, or inference~\cite{peherstorfer2018survey}. Full order forward models do not perform well in such many-query contexts, either due to excessive computational cost~(requiring the most powerful high performance computing architectures) or slow evaluation time (unacceptable in real-time contexts such as on-the-fly optimal control). In contrast to that of the \emph{big data} regime that dominates computer vision and other technological fields, only a relatively small amount of high resolution data can be generated from computer simulations or physical experiments in scientific applications. Fast approximate solvers built from this limited available data that can efficiently and accurately emulate the full order model would be highly advantageous. 

In this work, we demonstrate that the random feature model holds considerable potential for such a purpose. Resembling~\cite{lu2019deeponet,wu2020data} and the contemporaneous work in~\cite{bhattacharya2020pca, korolev2021two, li2020neural, o2020derivative}, we present a methodology for true function space learning of black-box input-output maps between a Banach space and separable Hilbert space. We formulate the approximation problem as supervised learning in infinite dimensions and show that the natural hypothesis space is a reproducing kernel Hilbert space associated with an operator-valued kernel. For a suitable loss functional, training the random feature model is equivalent to solving a finite-dimensional convex optimization problem. As a consequence of our careful construction of the method as mapping between Banach spaces, the resulting emulator naturally scales favorably with respect to the high input and output dimensions arising in practical, discretized applications; furthermore, it is shown to achieve small relative test error for two model problems arising from approximation of a semigroup and of the solution map corresponding to an elliptic PDE exhibiting parametric dependence on a coefficient function.

\subsection{Literature Review}
In recent years, two different lines of research have emerged that address PDE approximation problems with machine learning techniques. The first perspective takes a more traditional approach akin to point collocation methods from the field of numerical analysis. Here, the goal is to use a deep neural network (NN) to solve a prescribed initial boundary value problem with as high accuracy as possible. Given a point cloud in a spatio-temporal domain $ \tilde{D} $ as input data, the prevailing approach first directly parametrizes the PDE solution field as a NN and then optimizes the NN parameters by minimizing the PDE residual with respect to (w.r.t.) some loss functional~(see~\cite{raissi2019physics,sirignano2018dgm,weinan2018deep} and the references therein). To clarify, the object approximated with this novel method is a \emph{low-dimensional} input-output map $\tilde{D}\to\R $, i.e., the real-valued function that solves the PDE. This approach is mesh-free by definition but highly intrusive as it requires full knowledge of the specified PDE. Any change to the original formulation of the initial boundary value problem or related PDE problem parameters necessitates an (expensive) re-training of the NN solution. We do not explore this first approach any further in this article.

The second direction is arguably more ambitious: use a NN as an emulator for the infinite-dimensional mapping between an input parameter and the PDE solution itself or a functional of the solution, i.e., a quantity of interest; the latter is widely prevalent in uncertainty quantification problems. We emphasize that the object approximated in this setting, unlike in the aforementioned first approach, is an input-output map $ \cX\to\cY $, i.e., the PDE solution operator, where $ \cX,\, \cY $ are infinite-dimensional Banach spaces; this map is generally nonlinear. For an approximation-theoretic treatment of parametric PDEs in general, we refer the reader to the article of Cohen and DeVore~\cite{cohen2015approximation}. In applications, the solution operator is represented by a discretized forward model $ \R^{K} \to\R^{K}$, where $ K $ is the mesh size, and hence represents a \emph{high-dimensional} object. It is this second line of research that inspires our work.

Of course, there are many approaches to forward model reduction that do not explicitly involve machine learning ideas. The reduced basis method~(see \cite{barrault2004empirical,benner2017model,devore2014theoretical} and the references therein) is a classical idea based on constructing an empirical basis from data snapshots and solving a cheaper variational problem; it is still widely used in practice due to computationally efficient offline-online decompositions that eliminate dependence on the full order degrees of freedom. Recently, machine learning extensions to the reduced basis methodology, of both intrusive (e.g., projection-based reduced order models) and non-intrusive (e.g., model-free data only) type, have further improved the applicability of these methods~\cite{cheng2013data,gao2019non,hesthaven2018non,lee2020model,santo2019data}. However, the input-output maps considered in these works involve high dimension in only one of the input or output space, not both.
Other popular surrogate modeling techniques include Gaussian processes~\cite{williams2006gaussian}, polynomial chaos expansions~\cite{spanos1989stochastic}, and radial basis functions~\cite{wendland2004scattered}; yet, these are only practically suitable for problems with input space of low to moderate dimension. Classical numerical methods for PDEs may also represent the forward model $ \R^{K}\to\R^{K} $, albeit implicitly in the form a computer code~(e.g.: finite element, finite difference, finite volume methods). However, the approximation error is sensitive to $ K $ and repeated evaluations of this forward model often becomes cost prohibitive due to poor scaling with input dimension $ K $.

Instead, deep NNs have been identified as strong candidate surrogate models for parametric PDE problems due to their empirical ability to emulate high-dimensional nonlinear functions with minimal evaluation cost once trained. Early work in the use of NNs to learn the solution
operator, or vector field, defining ODEs and time-dependent PDEs,
may be found in the 1990s \cite{chen1995universal,Kev98,Kev92}. There are now more theoretical justifications for NNs breaking the \emph{curse of dimensionality}~\cite{korolev2021two,kutyniok2019theoretical,ma2019generalization}, leading to increased interest in PDE applications~\cite{adcock2020deep,geist2020numerical,opschoor2020deep,schwab2019deep}. A suite of work on data-driven discretizations of PDEs has surfaced that allow for identification of the governing model~\cite{bar2019learning,bigoni2020data,long2017pde,patel2020physics,stevens2020finitenet,trask2019gmls}; however, we note that only the operators appearing in the equation itself are approximated with these approaches, not the solution operator of the PDE. More in line with our focus in this article, architectures based on deep convolutional NNs have proven quite successful for learning elliptic PDE solution maps~(for example, see \cite{tripathy2018deep,winovich2019convpde,zhu2018bayesian}, which take an image-to-image regression approach). Other NNs have been used in similar elliptic problems for quantity of interest prediction~\cite{khoo2017solving}, error estimation~\cite{chen2020output}, or unsupervised learning~\cite{li2020variational}. Yet in all the approaches above, the architectures and resulting error are dependent on the mesh resolution. To circumvent this issue, the surrogate map must be well-defined on function space and independent of any finite-dimensional realization of the map that arises from discretization. This is not a new idea~(see~\cite{chen1995universal, rossi2005functional} or for functional data analysis,~\cite{kadri2016operator,micchelli2005learning}). The aforementioned reduced basis method is an example, as is the method of~\cite{chkifa2013sparse, cohen2015approximation}, which approximates the solution map with sparse Taylor polynomials and is proved to achieve optimal convergence rates in idealized settings. However, it is only recently that machine learning methods have been explicitly 
designed to operate in an infinite-dimensional setting, and there is little work in this direction~\cite{bhattacharya2020pca,li2020neural}. Here we propose the random feature model as another such method.

The random feature model~(RFM)~\cite{rahimi2008random,rahimi2008uniform,rahimi2008weighted}, detailed in~\Cref{sec:rfm}, is in some sense the simplest possible machine learning model; it may be viewed as an ensemble average of randomly parametrized functions: an expansion in a randomized basis. These \emph{random features} could be defined, for example, by randomizing the internal parameters of a NN. Compared to NN emulators with enormous learnable parameter counts~(e.g., $ O(10^{5}) $ to $ O(10^6) $, see~\cite{fan2020solving,feliu2020meta,li2020variational}) and methods that are intrusive or lead to nontrivial implementations~\cite{chkifa2013sparse,lee2020model,santo2019data}, the RFM is one of the simplest models to formulate and train (often $ O(10^{3})$ parameters, or fewer, suffice). The theory of the RFM for real-valued outputs is well developed, partly due to its close connection to kernel methods~\cite{bach2017equivalence,cao2019generalization,jacot2018neural,rahimi2008random,wendland2004scattered} and Gaussian processes~\cite{neal1996priors,williams1997computing}, and includes generalization rates and dimension-free estimates~\cite{ma2019generalization,rahimi2008uniform,sun2018random}. A quadrature viewpoint on the RFM provides further insight and leads to Monte Carlo sampling ideas~\cite{bach2017equivalence}; we remark on this further in~\Cref{sec:rfm}. As in modern deep learning practice, the RFM has also been shown to perform best when the model is over-parametrized~\cite{belkin2019reconciling}.  In a similar high-dimensional setting of relevance in this paper, the authors of~\cite{griebel2017reproducing,kempf2017kernel} theoretically investigated nonparametric kernel regression for parametric PDEs with real-valued solution map outputs. The specific random Fourier feature approach of Rahimi and Recht~\cite{rahimi2008random} was generalized in \cite{brault2016random} to the finite-dimensional matrix-valued kernel setting with vector-valued random Fourier features. However, most of these works require explicit knowledge of the kernel itself. Here our viewpoint is to work directly with random features as the basis for a standalone method, choosing them for their properties and
noting that they implicitly define a kernel, but not working directly with this kernel; furthermore, our work considers both infinite-dimensional input \emph{and} output spaces, not just one or the other. A key idea underlying our approach is to formulate the proposed random feature algorithm on infinite-dimensional space and only then discretize. This philosophy in algorithm development has been instructive in a number of areas in scientific computing, such as optimization \cite{hinze2008optimization} and the development of Markov chain Monte Carlo methodology \cite{cotter2013mcmc}. It has recently been promoted as a way of designing and analyzing algorithms within machine learning \cite{haber2017stable,ma2019machine,ruthotto2019deep,weinan2017proposal,weinan2019mean}, and our work may be understood within this general framework.

\subsection{Contributions}
Our primary contributions in this paper are now listed.
\begin{enumerate}
	\item We develop the random feature model, directly formulated on the function space level, for learning input-output maps between Banach spaces purely from data. As a method for parametric PDEs, the methodology is non-intrusive but also has the additional advantage that it may be used in settings where only data is available and no model is known.
	\item We show that our proposed method is more computationally tractable to both train and evaluate than standard kernel methods in infinite dimensions. Furthermore, we show that the method is equivalent to kernel ridge regression performed in a finite-dimensional space spanned by random features.
	\item We apply our methodology to learn the semigroup defined by the
	solution operator for viscous Burgers' equation and the coefficient-to-solution operator for the Darcy flow equation.
	\item We demonstrate, by means of numerical experiments, two mesh-independent approximation properties that are built into the proposed methodology: invariance of relative error to mesh resolution and evaluation ability on any mesh resolution.
\end{enumerate}

This paper is structured as follows. In~\Cref{sec:problem}, we communicate the mathematical framework required to work with the random feature model in infinite dimensions, identify an appropriate approximation space, and explain the training procedure. We introduce two instantiations of random feature maps that target physical science applications in~\Cref{sec:application} and detail the 
corresponding numerical results for these applications in~\Cref{sec:experiment}. We conclude in \Cref{sec:conclusion} with discussion and future work.

\section{Methodology}
\label{sec:problem}
In this work, the overarching problem of interest is the approximation of a map $ \Fd: \cX\to\cY $, where $ \cX,\,\cY $ are infinite-dimensional spaces of real-valued functions defined on some bounded open subset of $ \R^{d} $, and $\Fd $ is defined by $ a\mapsto \Fd(a)\defby u $, where $ u $ is the solution of a (possibly time-dependent) PDE and $a$ is an input function required to make the problem well-posed. Our proposed approach for this approximation, constructing a surrogate map $ F $ for the true map $ \Fd $, is 
data-driven, non-intrusive, and based on least squares. Least squares-based methods are integral to the random feature methodology as proposed in low dimensions~\cite{rahimi2008random,rahimi2008uniform} and generalized here to the infinite-dimensional setting; they have also been shown to work well in other algorithms for high-dimensional numerical approximation~\cite{beylkin2005algorithms,cohen2016optimal,doostan2009least}. Within the broader scope of reduced order modeling techniques~\cite{benner2017model}, the approach we adopt in this paper falls within the class of data-fit emulators. In its essence, our method interpolates the solution manifold
\begin{equation}\label{eqn:manifold}
\cM=\{u\in\cY: u=\Fd(a),\, a\in\cX\}\, .
\end{equation}
The solution map $ \Fd $, as the inverse of a differential operator, is often smoothing and admits a notion of compactness, i.e., the output space compactly embeds into the input space. Then, the idea is that $ \cM $ should have some compact, low-dimensional structure (intrinsic dimension). However, actually finding a model $ F $ that exploits this structure despite the high dimensionality of the truth map $ \Fd $ is quite difficult. Further, the effectiveness of many model reduction techniques, such as those based on the reduced basis method, are dependent on inherent properties of the map $ \Fd $ itself (e.g., analyticity), which in turn may influence the decay rate of the Kolmogorov width of the manifold $ \cM $~\cite{cohen2015approximation}. While such subtleties of approximation theory are crucial to
developing rigorous theory and provably convergent algorithms, we choose to work in the non-intrusive setting where knowledge of the map $ \Fd $ and its associated PDE are only obtained through measurement data, and hence detailed characterizations such as those aforementioned are essentially unavailable.

The remainder of this section introduces the mathematical preliminaries for our methodology. With the goal of operator approximation in mind, in~\Cref{sec:problem_form_SL} we formulate a supervised learning problem in an infinite-dimensional setting. We provide the necessary background on reproducing kernel Hilbert spaces in~\Cref{sec:rkhs} and then define the RFM in~\Cref{sec:rfm}. In~\Cref{sec:opt}, we describe the optimization principle which leads to algorithms for the RFM and an example problem in which $\cX$ and $\cY$ are one-dimensional vector spaces.

\subsection{Problem Formulation}
\label{sec:problem_form_SL}
Let $ \cX,\, \cY $ be real Banach spaces and $ \Fd:\cX\to\cY $ be a (possibly nonlinear) map. It is natural to frame the approximation of $ \Fd $ as a supervised learning problem. Suppose we are given training data in the form of input-output pairs $ \{a_i, y_i\}_{i=1}^{n}\subset \cX\times\cY$, where $ a_i\sim \nu$ i.i.d., $ \nu $ is a probability measure supported on $ \cX $, and $ y_i=\Fd(a_i) \sim \Fd_{\sharp}\nu$ with, potentially, noise added to the evaluations of $\Fd(\cdot)$. In the examples in this paper, the noise is viewed  as resulting from model error (the PDE does not perfectly represent the physics) or from discretization error (in approximating the PDE); situations in which the data acquisition process is inherently noisy can also be envisioned but are not studied here. We aim to build a parametric reconstruction of the true map $ \Fd $ from the data, that is, construct a model $ F:\cX\times\cP\to\cY $ and find $ \al^{\dagger}\in\cP\subseteq\R^{m} $ such that $ F(\cdot, \al^{\dagger})\approx \Fd $ are close as maps from $\cX$ to $\cY$ in some suitable sense. The natural number $ m $ here denotes the total number of model parameters. The standard approach to determine parameters in supervised learning is to first define a loss functional $ \ell:\cY\times\cY\to \R_{\geq 0} $ and then minimize the expected risk,
\begin{equation}\label{eqn:risk_expected}
\min_{\al\in\cP}\E^{a\sim\nu}\bigl[\ell\bigl(\Fd(a), F(a,\al)\bigr)\bigr]\, .
\end{equation}
With only the data $ \{a_i, y_i\}_{i=1}^{n} $ at our disposal, we approximate problem~\cref{eqn:risk_expected} by replacing $ \nu $ with the empirical measure $ \nu^{(n)}\defby\frac{1}{n}\sum_{j=1}^{n}\delta_{a_{j}} $, which leads to the empirical risk minimization problem
\begin{equation}\label{eqn:risk_empirical}
\min_{\al\in\cP}\dfrac{1}{n}\sum_{j=1}^{n}\ell\bigl(y_j, F(a_j,\al)\bigr)\, .
\end{equation}
The hope is that given minimizer $ \al^{(n)} $ of~\cref{eqn:risk_empirical} and $ \al^{\dagger} $ of~\cref{eqn:risk_expected}, $ F(\cdot,\al^{(n)}) $ well approximates $ F(\cdot, \al^{\dagger}) $, that is, the learned model \emph{generalizes} well; these ideas may be made rigorous with results from statistical learning theory~\cite{hastie2009elements}. Solving problem~\cref{eqn:risk_empirical} is called \emph{training} the model $ F $. Once trained, the model is then validated on a new set of i.i.d. input-output pairs previously unseen during the training process. This \emph{testing} phase indicates how well $ F $ approximates $ \Fd $. From here on out, we assume that $ (\cY, \ip{\cdot}{\cdot}_{\cY}, \norm{\cdot}_{\cY}) $ is a real separable Hilbert space and focus on the squared loss
\begin{equation}\label{eqn:loss_square}
\ell(y,y')\defby \dfrac{1}{2}\norm*{y-y'}^{2}_{\cY}\, .
\end{equation}
We stress that our entire formulation is in an infinite-dimensional setting and we will remain in this setting throughout the paper; as such, the random feature methodology we propose will inherit desirable discretization-invariant properties, to be observed in the numerical experiments of~\Cref{sec:experiment}.
\begin{notation}
	For a Borel measurable map $ G:\mathcal{U}\to\mathcal{V} $ between two Banach spaces $ \mathcal{U}$, $\mathcal{V} $ and a probability measure $ \pi $ supported on $ \mathcal{U} $, we denote the expectation of $ G $ under $ \pi $ by
	\begin{equation}\label{eqn:expect_nton}
	\E^{u\sim \pi}\bigl[G(u)\bigr]=\int_{\mathcal{U}} G(u)\pi(du)
	\end{equation}
	in the sense of Bochner integration (see, e.g., \cite{Dashti2017}, Sec. A.2). We will drop the domain of integration in situations where no confusion is caused by doing so.
\Endrmk\end{notation}

\subsection{Operator-Valued Reproducing Kernels}\label{sec:rkhs}
The random feature model is naturally formulated in a reproducing kernel Hilbert space (RKHS) setting, as our exposition will demonstrate in~\Cref{sec:rfm}. However, the usual RKHS theory is concerned with real-valued functions~\cite{aronszajn1950theory,berlinet2011reproducing,cucker2002mathematical,wendland2004scattered}. Our setting, with the output space $ \cY $ a separable Hilbert space, requires several ideas that generalize the real-valued case. We now outline these ideas with a review of operator-valued kernels; parts of the presentation that follow may be found in the references~\cite{bach2017equivalence,carmeli2006vector,micchelli2005learning}.

We first consider the special case $ \cY\defby \R $ for ease of exposition. A real RKHS is a Hilbert space $ (\cH, \ip{\cdot}{\cdot}_{\cH}, \norm{\cdot}_{\cH}) $ comprised of real-valued functions $ f:\cX\to\R $ such that the pointwise evaluation functional $ f\mapsto f(a) $ is bounded for every $ a\in\cX $. It then follows that there exists a unique, symmetric, positive definite kernel function $ k: \cX\times\cX\to\R $ such that for every $ a\in\cX $, $ k(\cdot, a)\in\cH $ and the \emph{reproducing kernel property} $ f(a)=\ip{k(\cdot,a)}{f}_{\cH} $ holds. These two properties are often taken as the definition of a RKHS. The converse direction is also true: every symmetric, positive definite kernel defines a unique RKHS~\cite{aronszajn1950theory}.

We now introduce the needed generalization of the reproducing property to the case of arbitrary real Hilbert spaces $ \cY $, as this result will motivate the construction of the RFM. Kernels in this setting are now operator-valued.
\begin{definition}\label{def:kernel_operatorvalued}
	Let $ \cX $ be a real Banach space and $ \cY $ a real separable Hilbert space. An \emph{\textbf{operator-valued kernel}} is a map
	\begin{equation}\label{eqn:kernel_operatorvalued}
	k:\cX\times \cX\to\cL(\cY,\cY)\, ,
	\end{equation}
	where $ \cL(\cY,\cY) $ denotes the Banach space of all bounded linear operators on $ \cY $, such that its adjoint satisfies $ k(a,a')^{*}=k(a',a) $ for all $ a,\, a'\in\cX $ and for every $ N\in\N $,
	\begin{equation}\label{eqn:kernel_psd}
	\sum_{i,j=1}^{N}\ip{y_i}{k(a_i,a_j)y_j}_{\cY}\geq 0
	\end{equation}
	for all pairs $ \{(a_i, y_i)\}_{i=1}^{N}\subset \cX\times\cY$.
\end{definition}

Paralleling the development for the real-valued case, an operator-valued kernel $ k $ also uniquely (up to isomorphism) determines an associated real RKHS $ \cH_{k}=\cH_{k}(\cX;\cY) $. Now, choosing a probability measure $ \nu $ supported on $ \cX $, we define a kernel integral operator $ T_k $ associated to $ k $ by
\begin{align}\label{eqn:integral_operator}
\begin{split}
T_{k}: L^{2}_{\nu}(\cX;\cY)&\to L^2_{\nu}(\cX;\cY)\\
F &\mapsto T_{k}F\defby \int k(\cdot,a')F(a')\nu(da')\, ,
\end{split}
\end{align}
which is non-negative, self-adjoint, and compact (provided $ k(a,a)\in\cL(\cY,\cY) $ is compact for all $ a\in\cX $~\cite{carmeli2006vector}). Let us further assume that all conditions needed  for $ T_{k}^{1/2} $ to be an isometry from $ L^{2}_{\nu} $ into $ \cH_{k} $ are satisfied, i.e., $ \cH_{k}=\im(T_k^{1/2}) $. Generalizing the standard Mercer theory~(see, e.g.,~\cite{bach2017equivalence,berlinet2011reproducing}), we may write the RKHS inner product as
\begin{equation}\label{eqn:ip_rkhs}
\ip{F}{G}_{\cH_{k}}=\ip{F}{T_{k}^{-1}G}_{L^2_{\nu}} \qfa F, G\in\cH_{k}\, .
\end{equation}
Note that while~\cref{eqn:ip_rkhs} appears to depend on the measure $ \nu $ on $ \cX $, the RKHS $ \cH_{k} $ is itself determined by the kernel without any reference to a measure~(see \cite{cucker2002mathematical}, Chp. 3, Thm. 4). With the inner product now explicit, we may directly deduce a reproducing property. A fully rigorous justification of the methodology is outside the scope of this article; however, we perform formal computations which provide intuition underpinning the methodology. To this end we fix $ a\in\cX $ and $ y\in\cY $. Then
\begin{align*}
\ip{k(\cdot, a)y}{T_k^{-1}F}_{L^{2}_{\nu}}
&=\int \ip*{k(a',a)y}{(T_{k}^{-1}F)(a')}_{\cY}\, \nu(da')\\
&=\int \ip*{y}{k(a,a')(T_{k}^{-1}F)(a')}_{\cY}\, \nu(da')\\
&=\ip[\Big]{y}{\int k(a,a')(T_{k}^{-1}F)(a')\, \nu(da')}_{\cY}\\
&=\ip{y}{F(a)}_{\cY}\, ,
\end{align*}
by using~\cref{def:kernel_operatorvalued} of operator-valued kernel and the fact that $ k(\cdot,a)y\in\cH_{k} $~(\cite{carmeli2006vector}). 
So, we deduce the following:
\begin{result}[Reproducing property for operator-valued kernels]
	Let $ F\in\cH_{k} $ be given. Then for every $ a\in\cX $ and $ y\in\cY $,
	\begin{equation}\label{eqn:rkprop_banach}
	\ip{y}{F(a)}_{\cY}=\ip{k(\cdot,a)y}{F}_{\cH_{k}}\, .
	\end{equation}
\end{result}

This identity, paired with a special choice of $ k $, is the basis of the random feature model in our abstract infinite-dimensional setting.

\subsection{Random Feature Model}\label{sec:rfm}
One could approach the approximation of target map $ \Fd:\cX\to\cY $ from the perspective of kernel methods. However, it is generally a difficult task to explicitly design operator-valued kernels of the form~\cref{eqn:kernel_operatorvalued} since the spaces $ \cX,\, \cY $ may be of different regularity, for example. Example constructions of operator-valued kernels studied in the literature include those taking value as diagonal operators, multiplication operators, or composition operators~\cite{kadri2016operator,micchelli2005learning}, but these all involve some simple generalization of scalar-valued kernels. Instead, the random feature model allows one to implicitly work with operator-valued kernels through the use of a \emph{random feature map} $ \varphi:\cX\times\Theta\to\cY $ and a probability measure $ \mu $ supported on Banach space $ \Theta $. The map $\varphi$ is assumed to be square 
integrable w.r.t. the product measure $\nu \times \mu$, i.e., $ \varphi\in L_{\nu\times\mu}^2(\cX\times\Theta;\cY) $, where $ \nu $ is the (sometimes a modeling choice at our discretion, sometimes unknown) data distribution on $ \cX $. Together, $ (\varphi, \mu) $ form a \emph{random feature pair}. With this setup in place, we now describe the connection
between random features and kernels; to this end, recall the following standard notation:
\begin{notation}
	Given a Hilbert space $ (H, \ip{\cdot}{\cdot}, \norm{\cdot}) $, the \emph{outer product} $ a\otimes b\in\cL(H,H) $ is defined by $ (a\otimes b)c=\ip{b}{c} a $ for any $ a,b,c\in H $.\Endrmk
\end{notation}

Given the pair $ (\varphi,\mu) $, consider maps $ k_{\mu}:\cX\times\cX\to\cL(\cY,\cY) $ of the form
\begin{equation}\label{eqn:kernel_expectation}
k_{\mu}(a,a')\defby \int \varphi(a;\theta)\otimes \varphi(a';\theta)\mu(d\theta)\, .
\end{equation}
Such representations need not be unique; different pairs $ (\varphi,\mu) $ may induce the same kernel $ k=k_{\mu} $ in \cref{eqn:kernel_expectation}. Since $ k_{\mu} $ may readily be shown to be an operator-valued kernel via~\cref{def:kernel_operatorvalued}, it defines a unique real RKHS $ \cH_{k_{\mu}}\subset L^{2}_{\nu}(\cX;\cY) $. Our approximation theory will be based on this space or finite-dimensional approximations thereof.
We now perform a purely formal but instructive calculation, following
from application of the reproducing property~\cref{eqn:rkprop_banach} 
to operator-valued kernels of the form~\cref{eqn:kernel_expectation}. 
Doing so leads to an integral representation of any $ F\in\cH_{k_{\mu}}$: 
for all $ a\in\cX,\, y\in\cY $,
\begin{align*}
\ip{y}{F(a)}_{\cY}=\ip{k_{\mu}(\cdot, a)y}{F}_{\cH_{k_{\mu}}}
&=\ip[\Big]{\int\ip{\varphi(a;\theta)}{y}_{\cY}\, \varphi(\cdot;\theta)\, \mu(d\theta)}{F}_{\cH_{k_{\mu}}}\\
&=\int \ip{\varphi(a;\theta)}{y}_{\cY}\ip{\varphi(\cdot;\theta)}{F}_{\cH_{k_{\mu}}}\mu(d\theta)\\
&=\int c_{F}(\theta)\ip{y}{\varphi(a;\theta)}_{\cY}\, \mu(d\theta)\\
&=\ip[\Big]{y}{\int c_{F}(\theta)\varphi(a;\theta)\mu(d\theta)}_{\cY}\, ,
\end{align*}
where the coefficient function $ c_{F}: \Theta\to\R $ is defined by
\begin{equation}\label{eqn:rfm_coeff_func}
c_{F}(\theta)\defby \ip{\varphi(\cdot;\theta)}{F}_{\cH_{k_{\mu}}}\, .
\end{equation}
Since $ \cY $ is Hilbert, the above holding for all $ y\in\cY$ 
implies the integral representation
\begin{equation}\label{eqn:rfm_integral_trans}
F = \int c_{F}(\theta)\varphi(\cdot;\theta)\mu(d\theta)\, .
\end{equation}
The formal expression \cref{eqn:rfm_coeff_func} for $c_F(\theta)$ needs careful interpretation (provided in \Cref{sec:integral_rkhs}). For instance, if $\varphi(\cdot;\theta)$ is a realization of a Gaussian process as in \Cref{ex:bb}, then $ \varphi(\cdot;\theta) \notin \cH_{k_{\mu}}$ with
probability one; indeed, in this case $c_F$ is defined only as an
$L^2_{\mu}$ limit. Nonetheless, the RKHS may be completely 
characterized by this integral representation. Define the map
\begin{align}\label{eqn:rf_integral_operator}
\begin{split}
\cA: L^{2}_{\mu}(\Theta;\R)&\to L^2_{\nu}(\cX;\cY)\\
c &\mapsto \cA c\defby \int c(\theta)\varphi(\cdot;\theta)\mu(d\theta) \, .
\end{split}
\end{align}
$ \cA $ may be shown to be a bounded linear operator that is a particular square root of $ T_{k_{\mu}} $ (\cref{sec:integral_rkhs}). We have the following result whose proof, provided in~\cref{sec:appendix_proofs}, is a straightforward
generalization of the real-valued case given in \cite{bach2017equivalence},
Sec. 2.2:

\begin{result}
\label{r:1}
	Under the assumption that $\varphi \in L^2_{\nu \times \mu}(\cX \times \Theta; \cY)$, the RKHS defined by the kernel $ k_{\mu} $ in~\cref{eqn:kernel_expectation} is precisely
	\begin{equation}\label{eqn:rkhs_image}
	\cH_{k_{\mu}}=\im(\cA)=\biggl\{\int c(\theta)\varphi(\cdot;\theta)\mu(d\theta): c\in L^{2}_{\mu}(\Theta;\R)\biggr\}\, .
	\end{equation} 
\end{result}

We stress that the integral representation of mappings in RKHS \cref{eqn:rkhs_image} is not unique since $ \cA $ is not injective in general. However, the particular choice $ c=c_{F} $ \cref{eqn:rfm_coeff_func} in representation \cref{eqn:rfm_integral_trans} does enjoy a sense of uniqueness as described in \Cref{sec:integral_rkhs}.

A central role in what follows is the approximation of measure $\mu$
by the empirical measure
\begin{equation}
\label{eqn:em}
\mu^{(m)} \defby  \frac{1}{m}\sum_{j=1}^{m}\delta_{\theta_{j}}\, ,\quad \theta_{j}\diid\mu\, .
\end{equation}
Given this, define $ k^{(m)}\defby k_{\mu^{(m)}} $ to be the empirical
approximation to $ k_{\mu} $:
\begin{equation}\label{eqn:kernel_empirical}
k^{(m)}(a,a')=\E^{\theta\sim\mu^{(m)}}\bigl[\varphi(a;\theta)\otimes\varphi(a';\theta)\bigr]=\dfrac{1}{m}\sum_{j=1}^{m}\varphi(a;\theta_{j})\otimes\varphi(a';\theta_{j})\, .
\end{equation}
Then we let $\cH_{k^{(m)}}$ be the unique RKHS induced by the 
kernel $k^{(m)}$; note that $ k^{(m)} $ and hence $ \cH_{k^{(m)}} $ are themselves random variables. The following characterization of $\cH_{k^{(m)}}$
is proved in~\Cref{sec:appendix_proofs}:

\begin{result}
\label{r:2}
Assume that $\varphi \in L^2_{\nu \times \mu}(\cX \times \Theta; \cY)$
and that the random features $\{\varphi(\cdot;\theta_j)\}_{j=1}^m$ are linearly independent in $L^2_{\nu}(\cX;\cY)$. Then, the RKHS $\cH_{k^{(m)}}$ is equal to the linear span of the $\{\varphi_j\defby \varphi(\cdot;\theta_j)\}_{j=1}^m$.
\end{result}

Applying a simple Monte Carlo sampling approach to elements in RKHS \cref{eqn:rkhs_image} by replacing probability measure $ \mu $ by empirical measure $\mu^{(m)}$ gives, for $ c\in L^2_{\mu} $,
\begin{equation}\label{eqn:MCRFM}
\frac{1}{m}\sum_{j=1}^{m}c(\theta_j)\varphi(\cdot;\theta_j)\approx \int c(\theta)\varphi(\cdot;\theta)\mu(d\theta)\, .
\end{equation}
This approximation achieves the Monte Carlo rate $ O(m^{-1/2})$ and, by virtue of \cref{r:2}, is in $\cH_{k^{(m)}}$. However, in the setting
of this work, the Monte Carlo approach does not give rise to
a practical method for learning a target map $ \Fd\in\cH_{k_{\mu}} $ because $\Fd$, $ k_{\mu} $, and $ \cH_{k_{\mu}} $ are all unknown; only the random feature pair $ (\varphi, \mu) $ is assumed to be given. Hence one cannot apply \cref{eqn:rfm_coeff_func} (or \cref{eqn:rkhs_norm_argmin}) to evaluate $ c=c_{\Fd} $ in \cref{eqn:MCRFM}. Furthermore, in realistic settings it may be that $ \Fd\not\in \cH_{k_{\mu}} $, which leads to an additional approximation gap not accounted for by the Monte Carlo method. To sidestep these  
difficulties, the RFM adopts a data-driven optimization approach 
to determine a different  approximation to $\Fd$, also from
the space $\cH_{k^{(m)}}$. We now define the RFM:

\begin{definition}\label{def:rfm}
Given probability spaces $(\cX,\mathcal{B}(\cX),\nu)$ and $(\Theta,\mathcal{B}(\Theta), \mu)$ with $\cX$, $ \Theta $ being real finite- or infinite-dimensional Banach spaces, 
real separable Hilbert space $\cY$, and 
$\varphi \in L^2_{\nu \times \mu}(\cX \times \Theta; \cY)$,
the \emph{\textbf{random feature model}} is the parametric map
\begin{align}\label{eqn:rfm}
\begin{split}
F_{m}: \cX\times\R^{m}&\to \cY\\
(a; \al)&\mapsto F_{m}(a; \al)\defby \dfrac{1}{m}\sum_{j=1}^{m}\al_{j}\varphi(a;\theta_{j})\, ,\quad \theta_{j}\diid\mu\, .
\end{split}
\end{align}
\end{definition}

We use the Borel $\sigma$-algebras $ \mathcal{B}(\cX)$ and $\mathcal{B}(\Theta) $ to define the probability spaces in the preceding definition. Our goal with the RFM is to choose parameters $\al\in\R^{m}$ so as to approximate mappings $ \Fd\in\cH_{k_{\mu}} $ (in the ideal setting) by mappings $ F_{m}(\cdot;\al)\in\cH_{k^{(m)}} $. The RFM is itself a random variable and may be viewed as a \emph{spectral method} since the randomized basis $ \varphi(\cdot;\theta) $ in the linear expansion~\cref{eqn:rfm} is defined on all of $ \cX$ $ \nu$-a.e.  
Determining the coefficient vector $\al$ from data obviates the difficulties associated with the Monte Carlo approach since the method only requires knowledge of the pair $ (\varphi,\mu) $ and knowledge of sample input-output pairs from target operator $\Fd$.

As written,~\Cref{eqn:rfm} is incredibly simple. It is clear that the choice of random feature map and measure pair $ (\varphi, \mu) $ will determine the quality of approximation. Most papers deploying these methods, including \cite{brault2016random,rahimi2008random,rahimi2008uniform}, take a kernel-oriented perspective by first choosing a kernel and then finding a random feature map to estimate this kernel. Our perspective, more aligned with \cite{rahimi2008weighted,sun2018random}, is the opposite in that we allow the choice of random feature map $ \varphi $ to implicitly \emph{define} the kernel via the formula~\cref{eqn:kernel_expectation} instead of picking the kernel first. This methodology also has implications for numerics: the kernel never explicitly appears in any computations, which leads to memory savings. It does, however, leave
open the question of characterizing the universality~\cite{sun2018random} of such kernels and the RKHS $\cH_{k_{\mu}}$ of mappings from $\cX$ to $\cY$ that underlies the approximation method; this is an important avenue for future work.

The close connection to kernels explains the origins of the RFM in the machine learning literature. Moreover, the RFM may also be interpreted in the context of neural networks \cite{neal1996priors,sun2018random,williams1997computing}. To see this explicitly, consider the setting where $\cX,\, \cY$ are both equal to the Euclidean space $\R$ and choose $ \varphi $ to be a family of hidden neurons $ \varphi_{\text{NN}}(a;\theta)\defby \sigma(\theta^{(1)}\cdot a + \theta^{(2)}) $. A single hidden layer NN would seek to find $\{(\alpha_j,\theta_j)\}_{j=1}^m$ in $\R \times \R^{2}$ so that
\begin{equation}
\label{eqn:vnn}
\frac{1}{m} \sum_{j=1}^m \alpha_j \varphi_{\text{NN}}(\cdot;\theta_j)
\end{equation}
matches the given training data $ \{a_i, y_i\}_{i=1}^{n}\subset \cX\times\cY$. More generally, and in arbitrary Euclidean spaces, one may allow
$ \varphi_{\text{NN}}(\cdot;\theta)$ to be any deep NN. 
However, while the RFM has the same \emph{form} as \cref{eqn:vnn}, there is a difference in the \emph{training}: the $ \theta_j $ are drawn i.i.d. from a probability measure and then fixed, and only the $\alpha_j$ are chosen to fit the training data. This connection is quite profound: given any deep NN with randomly initialized parameters $ \theta $, studies of the lazy training regime and neural tangent kernel \cite{cao2019generalization,jacot2018neural} suggest that adopting a RFM approach and optimizing over only $\al$ is quite natural, as it is observed that in this regime the internal NN parameters do not stray far from their random initialization during gradient descent whilst the last layer of parameters $\{\al_{j}\}_{j=1}^{m}$ adapt considerably. 

Once the feature parameters $ \{\theta_{j}\}_{j=1}^{m}$ are chosen at random and fixed, training the RFM $ F_{m} $ only requires optimizing over $ \al\in\R^{m} $ which, due to linearity of $F_m$ in $\al$, is a straightforward task to which we now turn our attention.

\subsection{Optimization}\label{sec:opt}
One of the most attractive characteristics of the RFM is its training procedure. With the $ L^2 $-type loss~\cref{eqn:loss_square} as in standard regression settings, optimizing the coefficients of the RFM with respect to the empirical risk \cref{eqn:risk_empirical} is a convex optimization problem, requiring only the solution of a finite-dimensional system of linear equations; the convexity also suggests the possibility of appending convex constraints (such as linear inequalities), although we do not pursue this here. Further, the kernels $ k_{\mu} $ or $ k^{(m)} $ are not required anywhere in the algorithm. We emphasize the simplicity of the underlying optimization tasks as they suggest the possibility of numerical implementation of the RFM into complicated black-box computer codes.

We now proceed to show that a regularized version of the optimization
problem~\cref{eqn:risk_empirical}--\cref{eqn:loss_square} arises naturally from approximation of a nonparametric regression problem defined over the RKHS $\cH_{k_{\mu}}.$ To this end, recall the supervised learning formulation in~\Cref{sec:problem_form_SL}. Given $ n $ i.i.d. input-output pairs $ \{a_i, y_i=\Fd(a_i)\}_{i=1}^{n}\subset \cX\times\cY$ as data, with the $a_i$ drawn from (possibly unknown)
probability measure $\nu$ on $ \cX $, the objective is to find an approximation $ \hat{F} $ to the map $ \Fd $. Let $ \cH_{k_{\mu}} $ be the hypothesis space and $ k_{\mu} $ its operator-valued reproducing kernel of the form~\cref{eqn:kernel_expectation}. The most straightforward learning algorithm in this RKHS setting is kernel ridge regression, also known as penalized least squares. This method produces a nonparametric model by finding a minimizer $ \hat{F} $ of
\begin{equation}\label{eqn:general_leastsquares}
\min_{F\in\cH_{k_{\mu}}}\biggl\{ \sum_{j=1}^{n}\dfrac{1}{2}\norm[\big]{y_j-F(a_j)}_{\cY}^{2}+\dfrac{\lambda}{2}\norm[\big]{F}_{\cH_{k_{\mu}}}^{2} \biggr\}\, ,
\end{equation}
where $ \lambda\geq 0 $ is a penalty parameter. By the representer theorem for operator-valued kernels (\cite{micchelli2005learning}, Theorems 2 and 4), the minimizer has the form
\begin{equation}\label{eqn:representer_thm}
\hat{F}=\sum_{j=1}^{n}k_{\mu}(\cdot,a_{j})\beta_{j}
\end{equation}
for some functions $\{\beta_{j}\}_{j=1}^{n}\subset\cY$. In practice, finding 
these $ n $ functions in the output space requires solving a block linear operator equation. For the high-dimensional PDE problems we consider in this work, solving such an equation may become prohibitively expensive from both operation count and memory required. A few workarounds were proposed in~\cite{kadri2016operator} such as certain diagonalizations, but these rely on simplifying assumptions that are somewhat limiting. More fundamentally, the representation of the solution in \cref{eqn:representer_thm} requires knowledge of the kernel $k_{\mu}$; in our setting we assume access only to the random feature pair $ (\varphi, \mu) $ which defines $k_{\mu}$ and not $k_{\mu}$ itself.

We thus explain how to make progress with this problem
given knowledge only of random features. Recall the empirical kernel given by \cref{eqn:kernel_empirical}, the RKHS $\cH_{k^{(m)}}$, and \cref{r:2}. The following result, proved in~\Cref{sec:appendix_proofs}, shows that a RFM hypothesis class with a penalized least squares empirical loss function in optimization problem \cref{eqn:risk_empirical}--\cref{eqn:loss_square} is equivalent to kernel ridge regression \cref{eqn:general_leastsquares} restricted to $\cH_{k^{(m)}}$. 
	
\begin{result}\label{res:minimizer_equiv}
Assume that $\varphi \in L^2_{\nu \times \mu}(\cX \times \Theta; \cY)$
and that the random features $\{\varphi(\cdot;\theta_j)\}_{j=1}^m$ are linearly independent in $L^2_{\nu}(\cX;\cY)$.
Fix $ \lambda\geq 0 $. Let $ \hat{\al}\in\R^{m} $ be the unique minimum norm 
solution of the following problem:  
\begin{equation}\label{eqn:opt_randomfeature}
\min_{\al\in\R^{m}}\biggl\{\sum_{j=1}^{n}\dfrac{1}{2}\norm[\bigg]{y_j-\dfrac{1}{m}\sum_{\ell=1}^{m}\al_{\ell}\varphi(a_j;\theta_{\ell})}_{\cY}^{2}+\dfrac{\lambda}{2m}\norm*{\al}_{2}^{2}\biggr\}\, .
\end{equation}
Then, the RFM defined by this choice $\al=\hat{\al}$ satisfies 
\begin{equation}\label{eqn:opt_equivalence}
F_{m}(\cdot;\hat{\al})=\argmin_{F\in\cH_{k^{(m)}}} \biggl\{\sum_{j=1}^{n}\dfrac{1}{2}\norm[\big]{y_j-F(a_j)}_{\cY}^{2}+\dfrac{\lambda}{2}\norm[\big]{F}_{\cH_{k^{(m)}}}^{2}\biggr\}\, .
\end{equation}
\end{result}

Solving the convex problem~\cref{eqn:opt_randomfeature} trains the RFM. The first order condition for a global minimizer leads to the normal equations
\begin{equation}\label{eqn:opt_normaleqn}
\dfrac{1}{m}\sum_{i=1}^{m}\sum_{j=1}^{n}\al_{i}\ip[\big]{\varphi(a_j;\theta_{i})}{\varphi(a_j;\theta_{\ell})}_{\cY} + \lambda \al_{\ell} =\sum_{j=1}^{n}\ip[\big]{y_j}{\varphi(a_j;\theta_{\ell})}_{\cY}
\end{equation}
for each $\ell\in\{1,\ldots,m\}$. This is an $ m $-by-$ m $ linear system of equations for $ \al\in\R^{m} $ that is standard to solve. In the case $ \lambda=0 $, the minimum norm solution may be written in terms of a pseudoinverse operator (see \cite{luenberger1997optimization}, Sec. 6.11).

\begin{example}[Brownian bridge]\label{ex:bb}
	We now provide a one-dimensional instantiation of the random feature model to illustrate the methodology. Take the input space as $ \cX\defby (0,1) $, output space $ \cY\defby \R $, input space measure $ \nu \defby U(0,1) $, and random parameter space $ \Theta\defby \R^{\infty} $. Denote the input by $ a=x\in\cX $. Then, consider the random feature map $ \varphi: (0,1)\times \R^{\infty}\to \R $ defined by the \emph{Brownian bridge}
	\begin{equation}\label{eqn:bb}
	\varphi(x;\theta)\defby \sum_{j\in\N}\theta^{(j)}(j\pi)^{-1}\sqrt{2}\sin(j\pi x)\, , \quad \theta^{(j)}\diid N(0,1)\, ,	
	\end{equation}
	where $ \theta\defby\{ \theta^{(j)} \}_{j\in\N} $ and $ \mu\defby N(0,1)\times N(0,1)\times\cdots$. For any realization of $ \theta \sim \mu$, the function $ \varphi(\cdot;\theta) $ is a Brownian motion constrained to zero at $ x=0 $ and $ x=1 $. The induced kernel $ k_{\mu}:(0,1)\times (0,1)\to\R $ is then simply the covariance function of this stochastic process:
	\begin{equation}\label{eqn:bb_kernel}
	k_{\mu}(x,x')=\E^{\theta\sim\mu}\bigl[\varphi(x;\theta)\varphi(x';\theta)\bigr]=\min\{x,x'\}-xx'\, .
	\end{equation}
	Note that $k_{\mu}$ is the Green's function for the negative Laplacian on $(0,1)$ with Dirichlet boundary conditions. Using this fact, we may explicitly characterize the associated RKHS $ \cH_{k_{\mu}} $ as follows. First, we have
	\begin{equation}\label{eqn:bb_covariance}
	T_{k_{\mu}}f = \int_{0}^{1}k_{\mu}(\cdot, y)f(y)\dd{y}=\Bigl(-\frac{d^2}{dx^2}\Bigr)^{-1}f\, ,
	\end{equation}
	where the the negative Laplacian has domain $H^2((0,1);\R)\cap H^1_0((0,1);\R) $. Viewing $ T_{k_{\mu}} $ as an operator from $ L^{2}((0,1);\R) $ into itself, from~\cref{eqn:ip_rkhs} we conclude, upon integration by parts, that
	\begin{equation}\label{eqn:bb_rkhsspace}
	\ip{f}{g}_{\cH_{k_{\mu}}}=\ip{f}{T_{k_{\mu}}^{-1}g}_{L^2}=\ip[\Big]{\frac{df}{dx}}{\frac{dg}{dx}}_{L^2}=\ip{f}{g}_{H_0^1} \qfa f, g\in\cH_{k_{\mu}}\, .
	\end{equation}
	Note that the last identity does indeed define an inner product on $H^1_0.$ By this formal argument we identify the RKHS $ \cH_{k_{\mu}} $ as the Sobolev space $ H_0^1((0,1);\R) $. Furthermore, Brownian bridge may be viewed as the Gaussian measure $N(0,T_{k_{\mu}})$.
	\begin{figure}[!htbp]
		\centering
		\begin{subfigure}[]{0.49\textwidth}
			\centering
			\includegraphics[width=\textwidth]{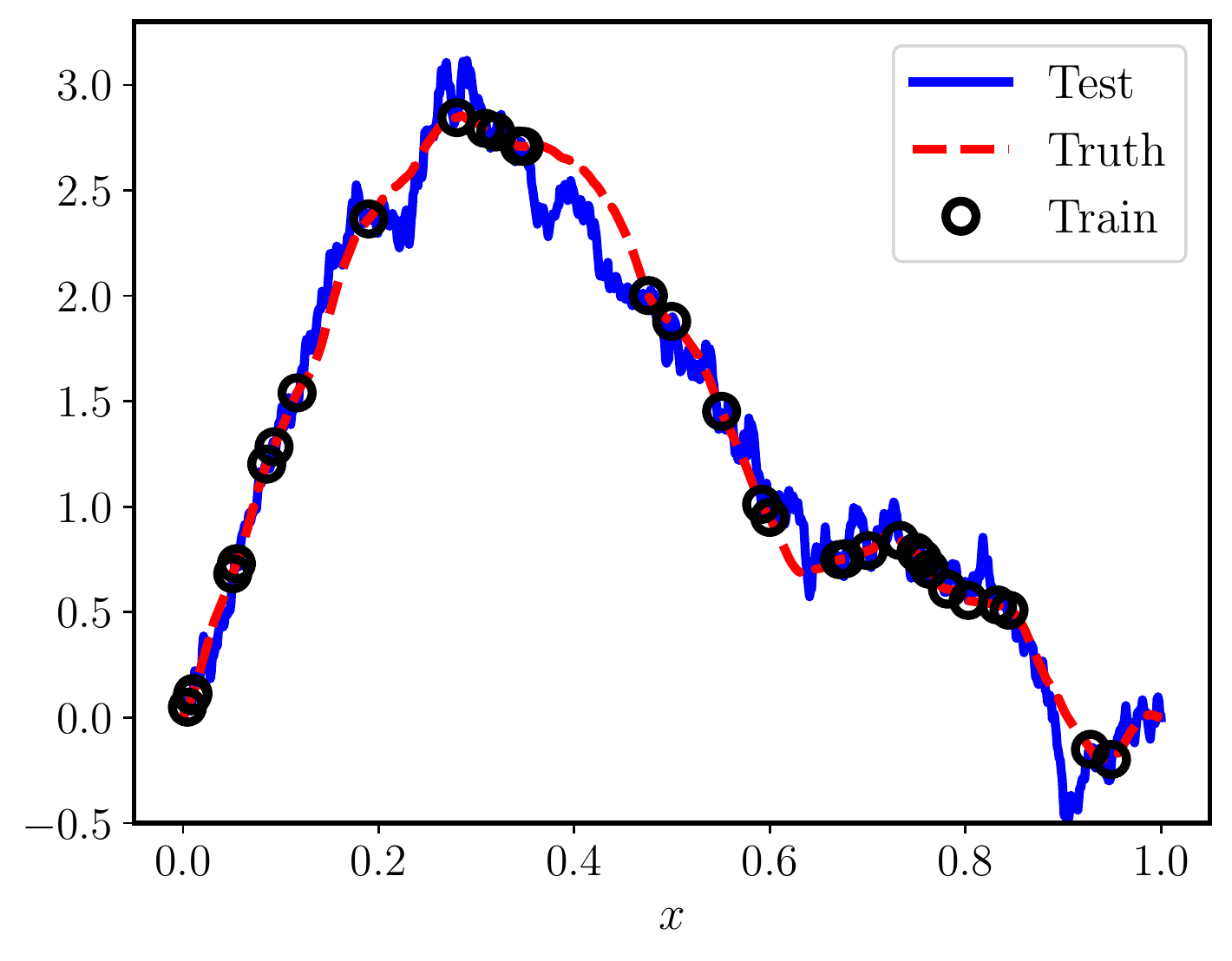}
			\caption{$ m=50 $}
			\label{fig:bb1}
		\end{subfigure}%
		\hfill%
		\begin{subfigure}[]{0.49\textwidth}
			\centering
			\includegraphics[width=\textwidth]{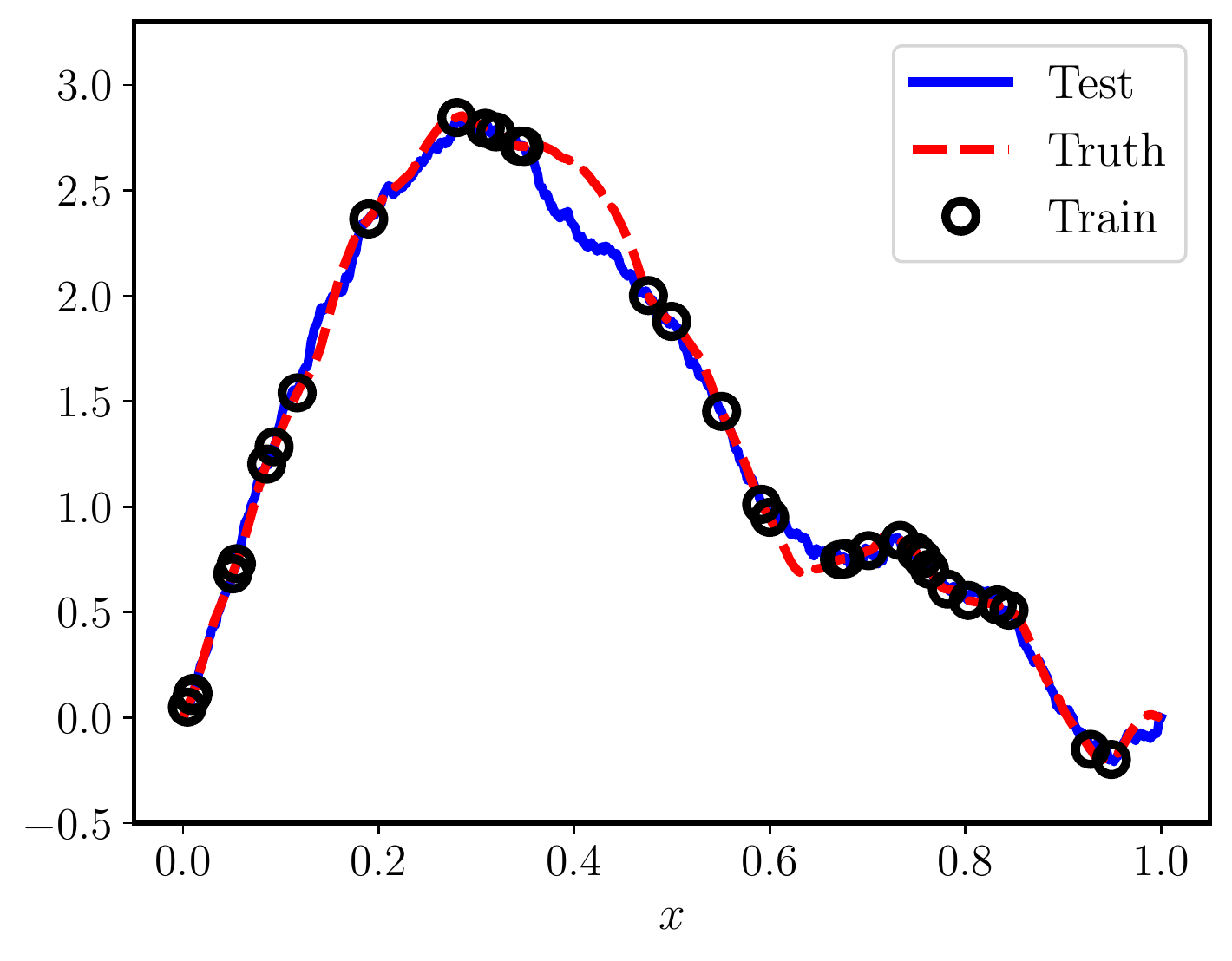}
			\caption{$ m=500 $}
			\label{fig:bb2}
		\end{subfigure}
		\begin{subfigure}[]{0.49\textwidth}
			\centering
			\includegraphics[width=\textwidth]{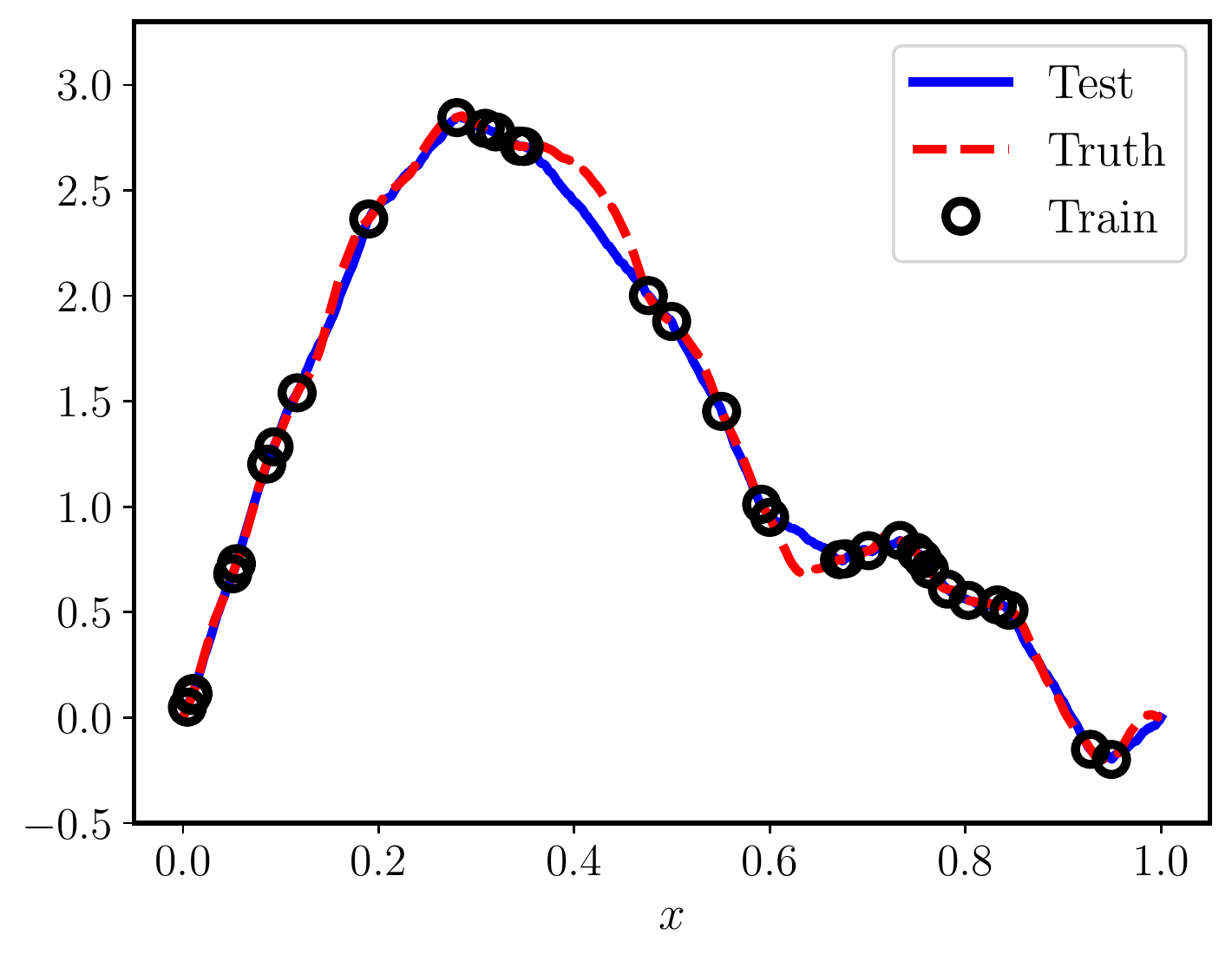}
			\caption{$ m=5000 $}
			\label{fig:bb3}
		\end{subfigure}%
		\hfill%
		\begin{subfigure}[]{0.49\textwidth}
			\centering
			\includegraphics[width=\textwidth]{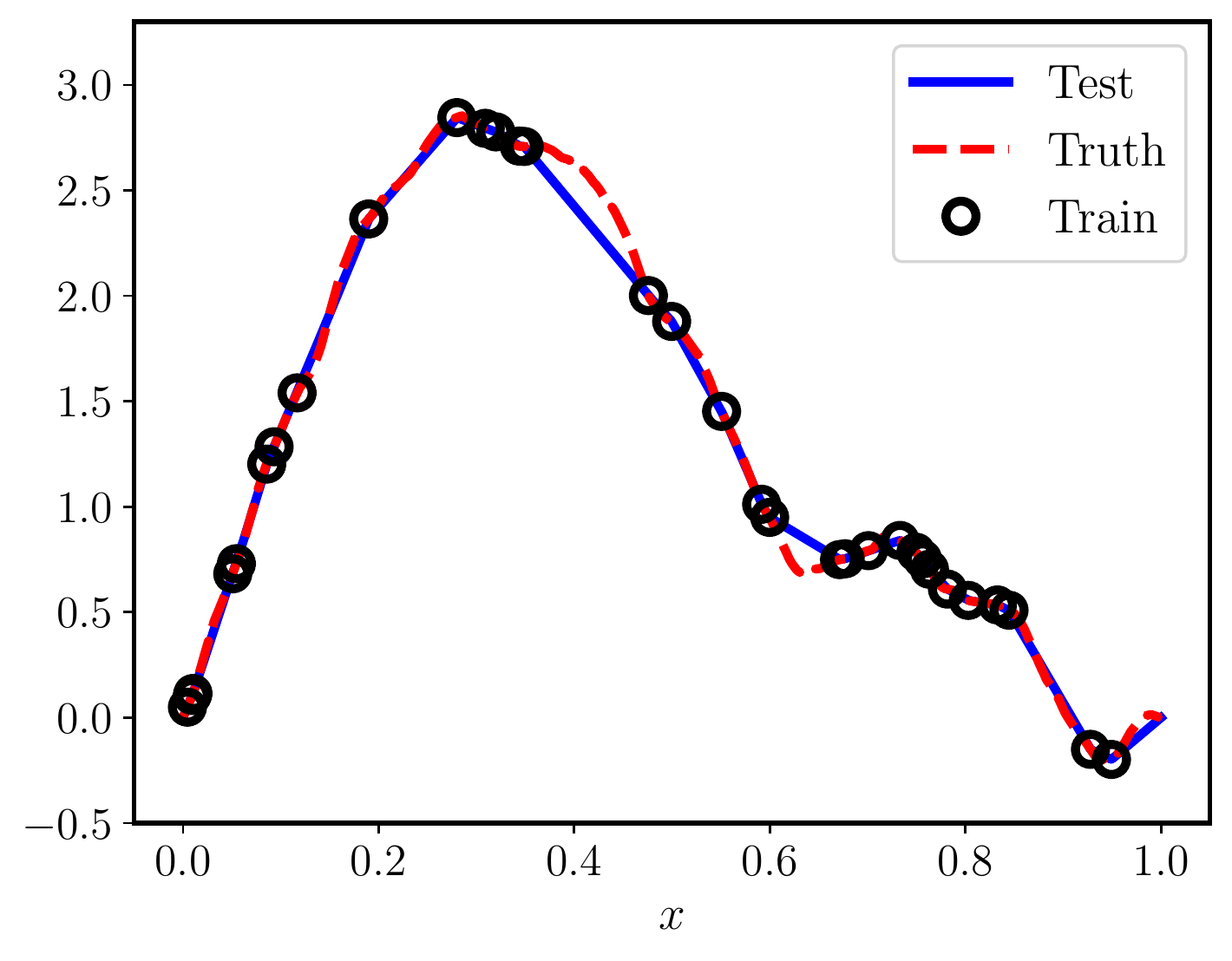}
			\caption{$ m=\infty $}
			\label{fig:bb4}
		\end{subfigure}
		\vspace{-5mm}
		\caption{Brownian bridge random feature model for one-dimensional input-output spaces with $ n=32 $ training points fixed and $ \lambda=0 $~(\Cref{ex:bb}): as $ m\to\infty $, the RFM approaches the nonparametric interpolant given by the representer theorem~(\Cref{fig:bb4}), which in this case is a piecewise linear approximation of the true function (an element of RKHS $ \cH_{k_{\mu}}=H_{0}^1 $, shown in red). Blue lines denote the trained model evaluated on test data points and black circles denote evaluation at training points.}
		\label{fig:bb_compare}
	\end{figure}
	Approximation using the RFM with the Brownian bridge random features is illustrated in~\Cref{fig:bb_compare}. Since $ k_{\mu}(\cdot,x) $ is a piecewise linear function, a kernel interpolation or regression method will produce a piecewise linear approximation. Indeed, the figure indicates that the RFM with $ n $ training points fixed approaches the optimal piecewise linear kernel interpolant as $ m\to\infty $ (see \cite{ma2019generalization} for a related theoretical result). 
	\Endrmk
\end{example}

The Brownian bridge \cref{ex:bb} illuminates a more fundamental idea. For this low-dimensional problem, an expansion in a deterministic Fourier sine basis would of course be more natural. But if we do not have a natural, computable orthonormal basis, then randomness provides a useful alternative representation; notice that the random features each include random combinations of the deterministic Fourier sine basis in this example. For the more complex problems that we study numerically in the next two sections, we lack knowledge of good, computable bases for general maps in infinite dimensions. The RFM approach exploits randomness to explore, implicitly discover the structure of, and represent, such maps. Thus we now turn away from this example of real-valued maps defined on a subset of the real line and instead consider the use of random features to represent maps between spaces of functions.

\section{Application to PDE Solution Maps}
\label{sec:application}
In this section, we design the random feature maps $ \varphi: \cX\times\Theta\to\cY $ and measures $ \mu $ for the RFM approximation of two particular PDE parameter-to-solution maps: the evolution semigroup of viscous Burgers' equation in \Cref{sec:burg_formulation} and the coefficient-to-solution operator for the Darcy problem in \Cref{sec:darcy_formulation}. It is well known to kernel method practitioners that the choice of kernel (which in this work follows from the choice of $ (\varphi,\mu) $) plays a central role in the quality of the function reconstruction. While our method is purely data-driven and requires no knowledge of the governing PDE, we take the view that any prior knowledge can, and should, be introduced into the design of $ (\varphi, \mu)$. However, the question of how to automatically determine good random feature pairs for a particular problem or dataset, inducing data-adapted kernels, is open. The maps $ \varphi$ that we choose to employ are nonlinear in both arguments. We also detail the probability measure $ \nu $ on the input space $ \cX $ for each of the two PDE applications; this choice is crucial because while we desire the trained RFM to transfer to arbitrary out-of-distribution inputs from $ \cX $, we can in general only expect the learned map to perform well when restricted to inputs statistically similar to those sampled from $ \nu $.

\subsection{Burgers' Equation: Formulation}\label{sec:burg_formulation}
Viscous Burgers' equation in one spatial dimension is representative of the advection-dominated PDE problem class in some regimes; these time-dependent equations are not conservation laws due to the presence of small dissipative terms, but nonlinear transport still plays a central role in the evolution of solutions. The initial value problem we consider is
\begin{equation}\label{eqn:burgers_ibvp}
\begin{cases}
\begin{alignedat}{2}
\frac{\partial u}{\partial t}+\frac{\partial}{\partial x}\left(\frac{u^2}{2}\right)-\ep \frac{\partial^2 u}{\partial x^2}&=f\, \ \ &&\text{in } (0,\infty)\times(0,1)\,,\\
u(\cdot, 0)=u(\cdot, 1)\, ,\quad \frac{\partial u}{\partial x}(\cdot, 0)&=\frac{\partial u}{\partial x}(\cdot, 1)\, \ \  &&\text{in } (0,\infty)\,,\\
u(0,\cdot)&=a\, \ \ &&\text{in } (0,1)\, ,
\end{alignedat}
\end{cases}
\end{equation}
where $ \ep>0 $ is the viscosity (i.e., diffusion coefficient) and we have imposed periodic boundary conditions. The initial condition $ a $ serves as the input and is drawn according to a Gaussian measure defined by
\begin{equation}\label{eqn:prior_burgers}
a\sim \nu\defby N(0,C)
\end{equation}
with Mat\'ern-like covariance operator \cite{dunlop2017hierarchical,matern2013spatial}
\begin{equation}\label{eqn:prior_covariance}
C\defby\tau^{2\al-d}(-\lap +\tau^{2}\id)^{-\al}\, ,
\end{equation}
where $d=1$ and the negative Laplacian $ -\lap $ is defined over $\mathbb{T}^1=[0,1]_{\mathrm{per}}$ and restricted to functions which integrate to zero over $ \mathbb{T}^1 $. The hyperparameter $ \tau\geq 0 $ is an inverse length scale  and $ \al>1/2 $ controls the regularity of the draw. Such $ a $ are almost surely H\"older and Sobolev regular with exponent up to $ \al - 1/2 $ (\cite{Dashti2017}, Thm. 12, pg. 338), so in particular $ a\in \cX\defby L^{2}(\mathbb{T}^1;\R) $. Then for all $ \ep>0 $, the unique global solution $ u(t,\cdot) $ to~\cref{eqn:burgers_ibvp} is real analytic for all $ t>0 $~(see \cite{kiselev2008blow}, Thm. 1.1). Hence, setting the output space to be $ \cY\defby  H^{s}(\mathbb{T}^1;\R) $ for any $ s>0 $, we may define the solution map

\begin{align}\label{eqn:solnmap_burg}
\begin{split}
\Fd: L^2 &\to H^{s} \\
a&\mapsto \Fd(a)\defby \Psi_{T}(a)=u(T,\cdot)\, ,
\end{split}
\end{align}
where $ \{\Psi_{t}\}_{t>0} $ forms the solution operator semigroup for~\cref{eqn:burgers_ibvp} and we fix the final time $ t=T>0 $. The map $ \Fd $ is smoothing and nonlinear.

\begin{figure}[!htbp]
	\centering
	\begin{subfigure}[]{0.49\textwidth}
		\centering
		\includegraphics[width=\textwidth]{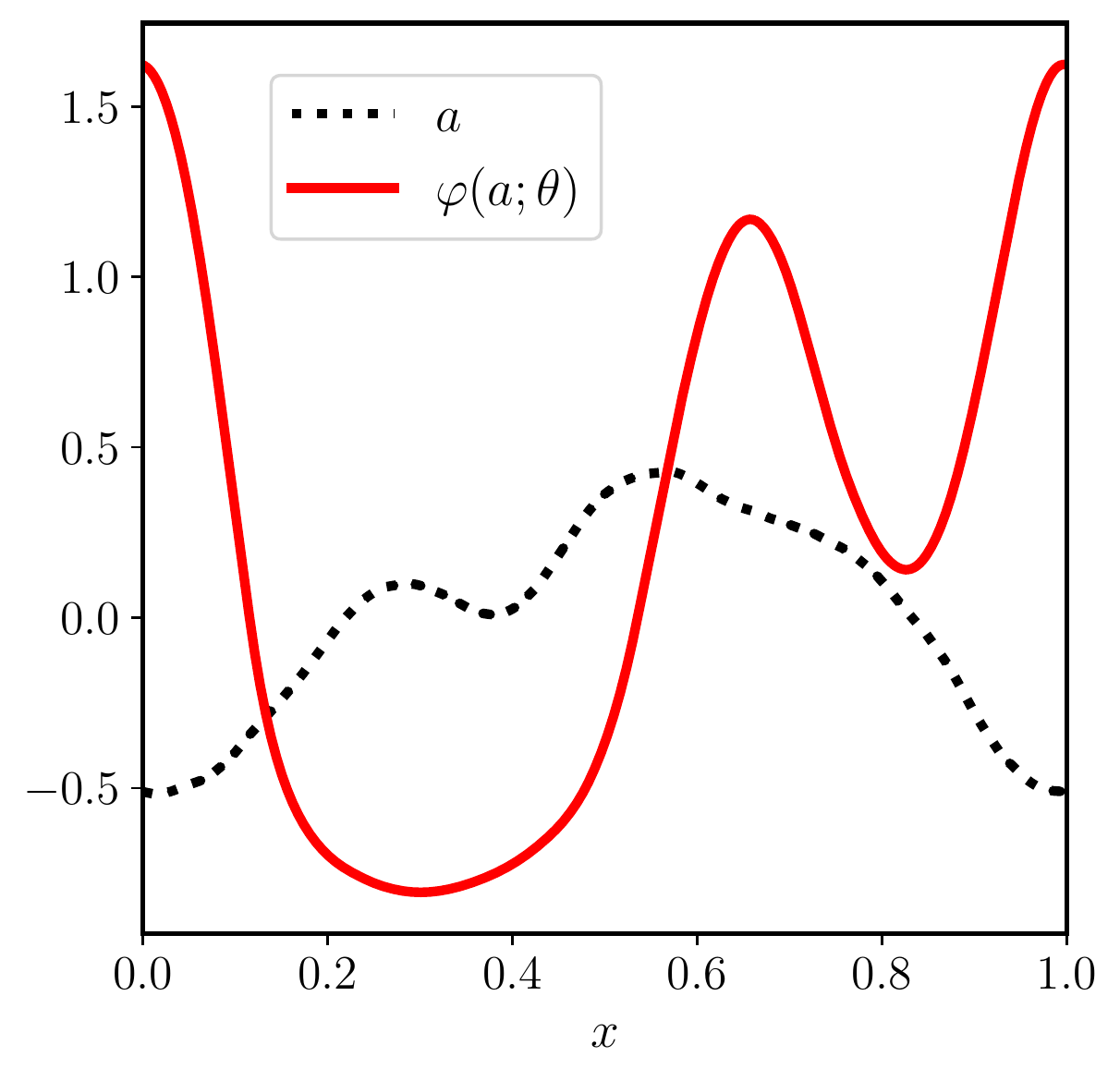}
		\caption{}
		\label{fig:rf_sample_burg}
	\end{subfigure}%
	\hfill%
	\begin{subfigure}[]{0.49\textwidth}
		\centering
		\includegraphics[width=\textwidth]{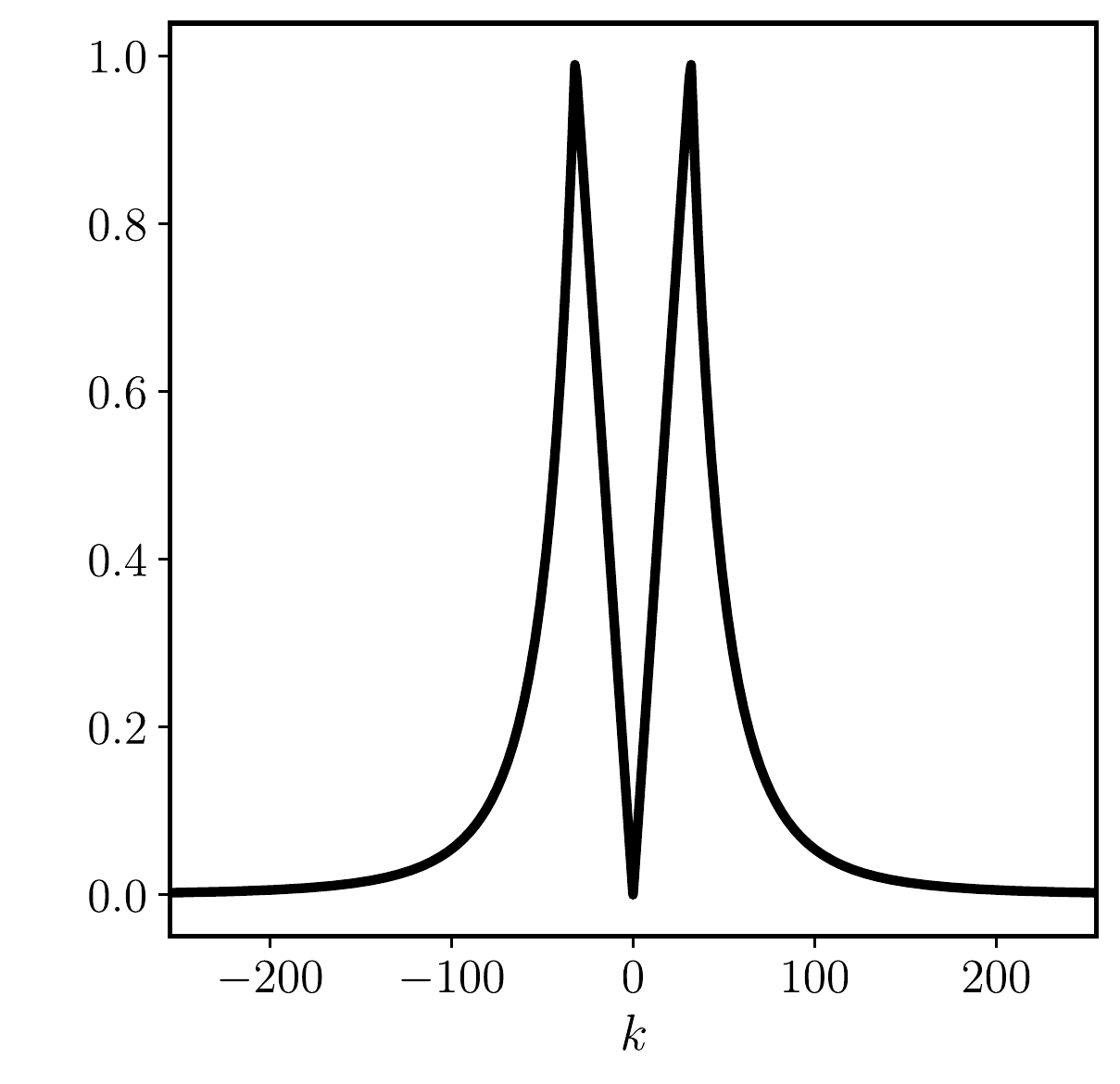}
		\caption{}
		\label{fig:filter_func1}
	\end{subfigure}
	\vspace{-5mm}
	\caption{Random feature map construction for Burgers' equation: \Cref{fig:rf_sample_burg} displays a representative input-output pair for the random feature $ \varphi(\cdot;\theta) $, $ \theta\sim\mu $~\cref{eqn:rf_fourier}, while \Cref{fig:filter_func1} shows the filter $k\mapsto \chi(k) $ for $ \delta=0.0025 $ and $ \beta=4 $~\cref{eqn:filter}.}
	\label{fig:rf_and_filter}
\end{figure}

We now describe a random feature map for use in the RFM~\cref{eqn:rfm} that we call \emph{Fourier space random features}. Let $ \cF $ denote the Fourier transform over spatial domain $ \mathbb{T}^1 $ and define $ \varphi: \cX\times\Theta\to\cY $ by 
\begin{equation}\label{eqn:rf_fourier}
\varphi(a;\theta)\defby\sigma\left(\cF^{-1}(\chi\cF a\cF \theta)\right)\, ,
\end{equation}
where $\sigma(\cdot)$, the $ \ON{ELU} $ function defined below, is defined as a mapping on $\mathbb{R}$ and applied pointwise to functions. Viewing $ \Theta\subseteq L^{2}(\mathbb{T}^1;\R) $, the randomness enters through $ \theta\sim \mu\defby N(0,C') $ with $ C' $ the same covariance operator as in~\cref{eqn:prior_covariance} but with potentially different inverse length scale and regularity, and the \emph{wavenumber filter function} $ \chi:\Z\to\R_{\geq 0} $ is
\begin{equation}\label{eqn:filter}
\chi(k)\defby \sigma_{\chi}(2\pi\abs{k}\delta)\, , \quad \sigma_{\chi}(r)\defby \max\bigl\{0, \min\{2r, (r+1/2)^{-\beta}\}\bigr\}\, ,
\end{equation}
where $ \delta,\, \beta>0 $. The map $ \varphi(\cdot;\theta) $ essentially performs a filtered random convolution with the initial condition. \Cref{fig:rf_sample_burg} illustrates a sample input and output from $ \varphi $. Although simply hand-tuned for performance and not optimized, the filter $ \chi $ is designed to shuffle energy in low to medium wavenumbers and cut off high wavenumbers~(see~\cref{fig:filter_func1}), reflecting our prior knowledge of solutions to~\cref{eqn:burgers_ibvp}. 

We choose the activation function $ \sigma $ in~\cref{eqn:rf_fourier} to be the exponential linear unit
\begin{equation}\label{eqn:activation_elu}
\ON{ELU}(r)\defby 
\begin{cases}
\begin{alignedat}{2}
r&\, , \ \ \ &&r\geq 0\\
e^{r}-1&\, , \ \ \ &&r<0\, .
\end{alignedat}
\end{cases}
\end{equation}
$ \ON{ELU} $ has successfully been used as activation in other  machine learning frameworks for related nonlinear PDE problems~\cite{lee2020model, patel2018nonlinear,patel2020physics}. We also find $ \ON{ELU} $ to perform better in the RFM framework over several other choices including $ \ON{ReLU}(\cdot)$, $\tanh(\cdot)$, $\ON{sigmoid}(\cdot)$, $\sin(\cdot) $, $ \ON{SELU}(\cdot) $, and $ \ON{softplus}(\cdot) $. Note that the pointwise evaluation of $ \ON{ELU} $ in~\cref{eqn:rf_fourier} will be well defined, by Sobolev embedding, for $ s>1/2 $ sufficiently large in the definition of $ \cY=H^{s} $. Since the solution operator maps into $H^s$ for any $s>0$, this does not constrain the method.

\subsection{Darcy Flow: Formulation}\label{sec:darcy_formulation}
\begin{figure}[!htbp]
	\centering
	\begin{subfigure}[]{0.49\textwidth}
		\centering
		\includegraphics[width=\textwidth]{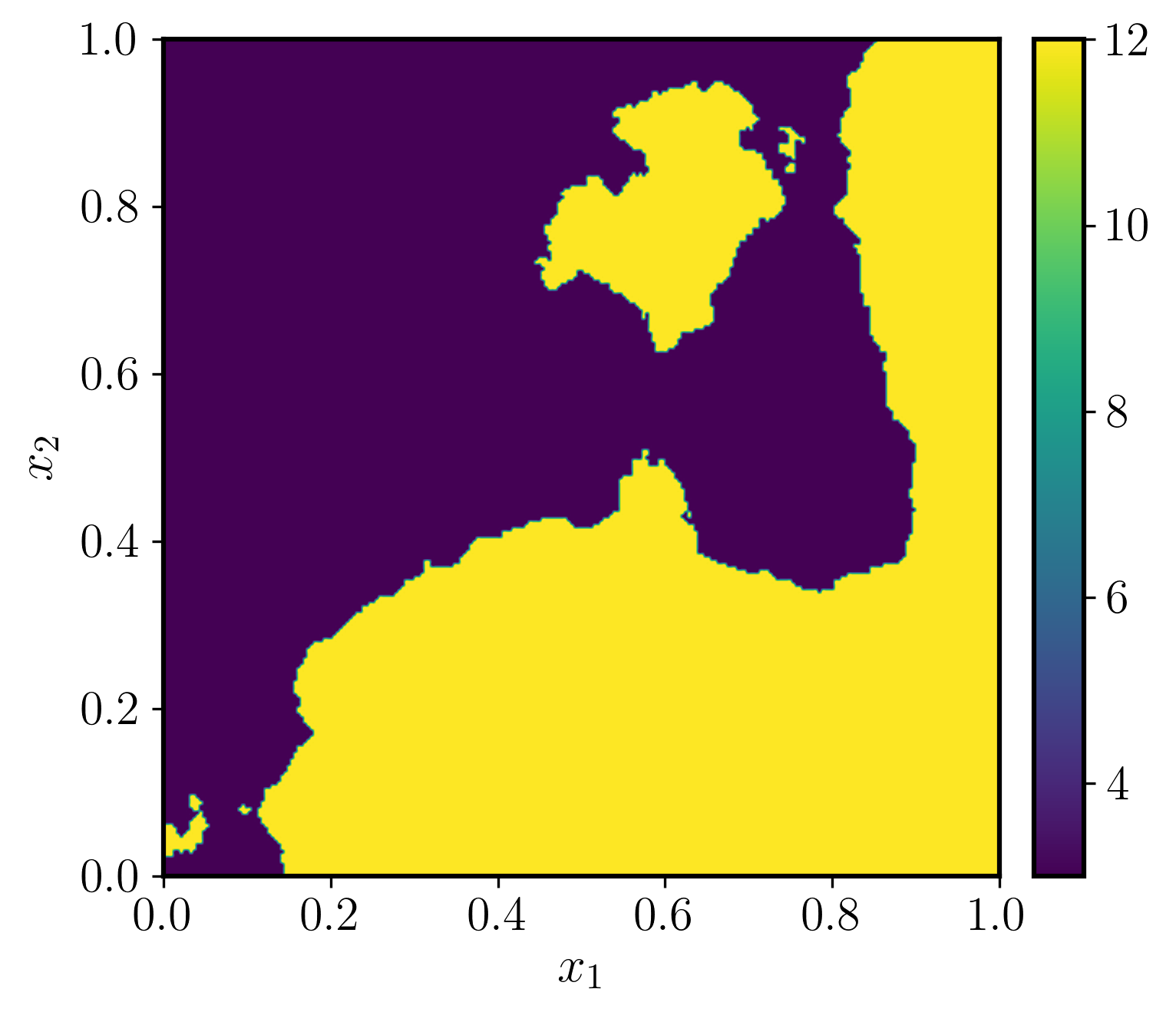}
		\caption{$ a\sim \nu $}
		\label{fig:coef_darcy}
	\end{subfigure}%
	\hfill%
	\begin{subfigure}[]{0.49\textwidth}
		\centering
		\includegraphics[width=\textwidth]{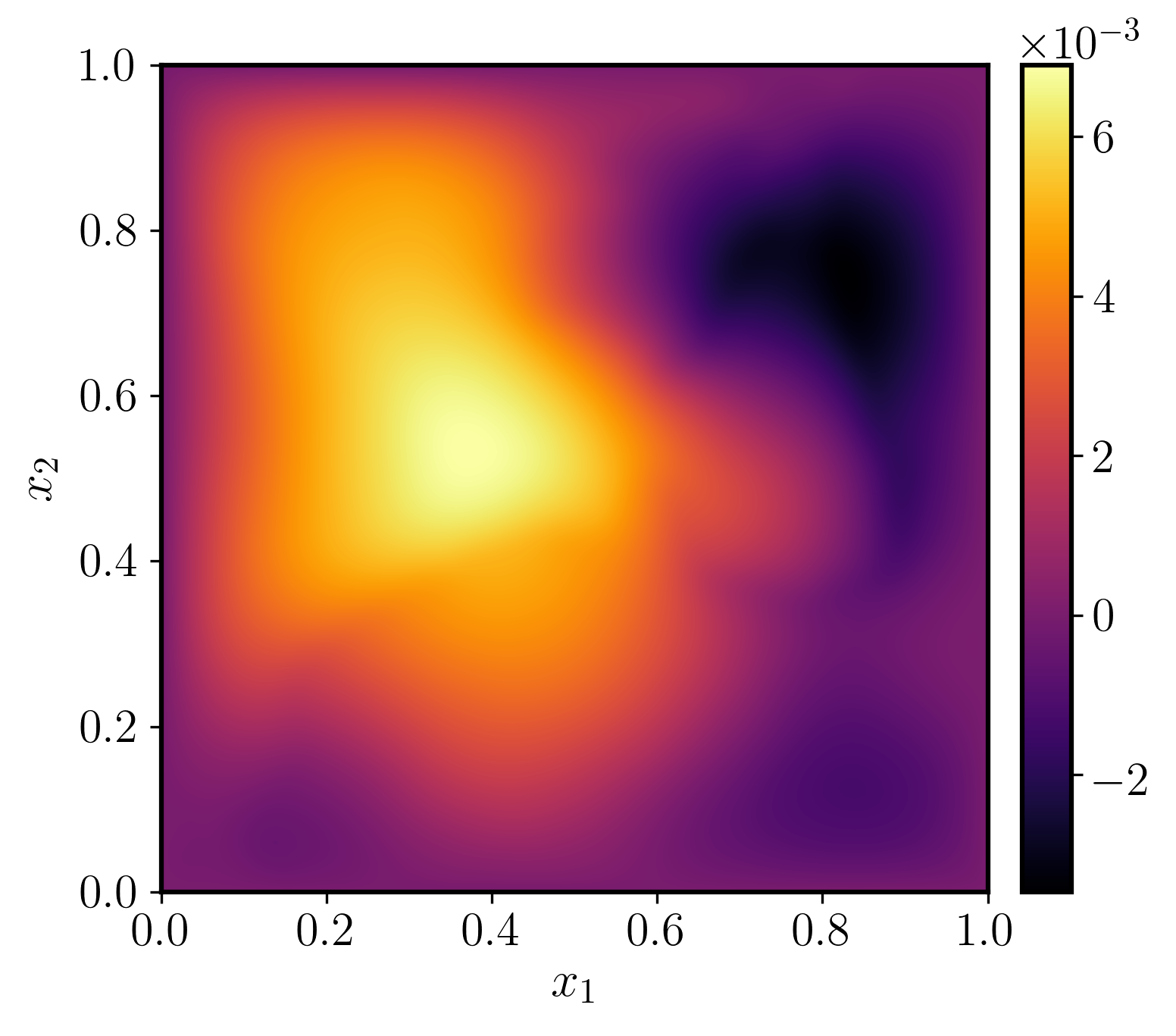}
		\caption{$ \varphi(a;\theta)\, , \ \theta\sim\mu$}
		\label{fig:rf_darcy}
	\end{subfigure}
	\vspace{-5mm}
	\caption{Random feature map construction for Darcy flow: \Cref{fig:coef_darcy} displays a representative input draw $ a $ with $ \tau=3,\, \al=2 $  and $ a^{+}=12,\, a^{-}=3 $; \Cref{fig:rf_darcy} shows the output random feature $ \varphi(a;\theta) $~(\Cref{eqn:rf_predictor_corrector}) taking the coefficient $ a $ as input. Here, $ f\equiv 1 $, $ \tau'=7.5,\, \al'=2 $, $ s^{+}=1/a^{+} $, $ s^{-}=-1/a^{-} $, and $ \delta = 0.15 $.}
	\label{fig:rf_and_coef_darcy}
\end{figure}

Divergence form elliptic equations~\cite{gilbarg2015elliptic} arise in a variety of applications, in particular, the groundwater flow in a porous medium governed by Darcy's law~\cite{bear2012fundamentals}. This linear elliptic boundary value problem reads
\begin{equation}\label{eqn:darcy_flow}
\begin{cases}
\begin{alignedat}{2}
-\nabla\cdot(a\nabla u)&=f\,  \ \ &&\text{in } D\,,\\
u&=0\,  \ \ &&\text{on } \partial D\, ,
\end{alignedat}
\end{cases}
\end{equation}
where $ D $ is a bounded open subset in $ \R^{d} $, $f$ represents sources
and sinks of fluid, $a$ the permeability of the porous medium, and $u$ the piezometric head; all three functions map $D$ into $\R$ and, in addition, $a$ is strictly positive almost everywhere in $D$. We work in a setting where $f$ is fixed and consider the input-output map defined by $a \mapsto u$. The measure $\nu$ on $a$ is a high contrast level set prior 
constructed as the pushforward of a Gaussian measure:
\begin{equation}\label{eqn:prior_levelset}
a\sim \nu\defby \psi_{\sharp}N(0,C)\, .
\end{equation}
Here $ \psi:\R\to\R $ is a threshold function defined by
\begin{equation}\label{eqn:prior_function_push}
\psi(r)\defby a^{+}\one_{(0,\infty)}(r)+a^{-}\one_{(-\infty,0)}(r)\, ,\quad 0<a^{-}\leq a^{+}<\infty\, ,
\end{equation}
applied pointwise to functions, and the covariance operator $ C $ is given in~\cref{eqn:prior_covariance} with $ d=2 $ and homogeneous Neumann boundary conditions on $ -\lap $. That is, the resulting coefficient $ a $ almost surely takes only two values ($ a^{+} $ or $ a^{-} $) and, as the zero level set of a Gaussian random field, exhibits random geometry in the physical domain $ D $. It follows that $ a\in L^{\infty}(D;\R_{\geq 0}) $ almost surely. Further, the size of the contrast ratio $ a^{+}/a^{-} $ measures the scale separation of this elliptic problem and hence controls the difficulty of reconstruction~\cite{bernardi2000adaptive}. See~\Cref{fig:coef_darcy} for a representative sample.

Given $ f\in L^{2}(D;\R) $, the standard Lax-Milgram theory may be applied to show that for coefficient $a\in\cX\defby L^{\infty}(D;\R_{\geq 0})$, there exists a unique weak solution $ u\in \cY\defby H_{0}^{1}(D;\R) $ for~\Cref{eqn:darcy_flow}~(see, e.g.,~Evans~\cite{evans2010partial}). Thus, we define the ground truth solution map 
\begin{align}\label{eqn:solnmap_darcy}
\begin{split}
\Fd: L^{\infty} &\to H_{0}^{1} \\
a&\mapsto \Fd(a)\defby u\, .
\end{split}
\end{align}
Although the PDE~\cref{eqn:darcy_flow} is linear, the solution map $ \Fd $ is nonlinear.

We now describe the chosen random feature map for this problem, which we call \emph{predictor-corrector random features}. Define $ \varphi: \cX\times\Theta\to\cY$ by $ \varphi(a;\theta)\defby p_1 $ such that
\begin{subequations}\label{eqn:rf_predictor_corrector}
	\begin{align}
	-\lap p_0&=\dfrac{f}{a}+\sigma_{\gamma}(\theta_1)\,,\label{eqn:predictor}\\
	-\lap p_1&=\dfrac{f}{a}+\sigma_{\gamma}(\theta_2)+\nabla(\log a)\cdot\nabla p_0\, ,\label{eqn:corrector}
	\end{align}
\end{subequations}
where the boundary conditions are homogeneous Dirichlet, $ \theta=(\theta_1,\theta_2) \sim \mu\defby\mu'\times\mu'$ are two Gaussian random fields each drawn from $ \mu'\defby N(0, C') $, $ f $ is the source term in~\cref{eqn:darcy_flow}, and $ \gamma =(s^{+}, s^{-}, \delta)$ are parameters for a thresholded sigmoid $ \sigma_{\gamma}:\R\to\R $,
\begin{equation}\label{eqn:sigmoidal_func_defn}
\sigma_{\gamma}(r)\defby \dfrac{s^{+}-s^{-}}{1+e^{-r/\delta}}+s^{-}\, ,
\end{equation}
and extended as a Nemytskii operator when applied to $\theta_1(\cdot)$ 
or $\theta_2(\cdot)$. We view $ \Theta\subseteq L^2(D;\R)\times L^2(D;\R) $. In practice, since $ \nabla a $ is not well-defined when drawn from the level set measure, we replace $ a $ with $a_{\ep} $, where $ a_{\ep} \defby v(1) $ is a smoothed version of $ a $ obtained by evolving the following linear heat equation for one time unit:
\begin{equation}\label{eqn:heat_mollify_neumann}
\begin{cases}
\begin{alignedat}{2}
\dv{v}{t}&=\eta\lap v \ \ &&\text{in } (0,1)\times D\,,\\
n\cdot \nabla v &= 0 \ \ &&\text{on } (0,1)\times\partial D\,,\\
v(0) &= a \ \ &&\text{in } D\,,
\end{alignedat}
\end{cases}
\end{equation}
where $ n $ is the outward unit normal vector to $ \partial D $. An example of the response $ \varphi(a;\theta) $ to a piecewise constant input $ a\sim\nu$ is shown in~\Cref{fig:rf_and_coef_darcy} for some $ \theta\sim\mu $.

We remark that by removing the two random terms involving $ \theta_1,\, \theta_2 $ in~\cref{eqn:rf_predictor_corrector}, we obtain a remarkably accurate surrogate model for the PDE. This observation is representative of a more general iterative method, a predictor-corrector type iteration, for solving the Darcy equation~\cref{eqn:darcy_flow}, whose convergence depends on the size of $ a $. The map $ \varphi $ is essentially a random perturbation of a single step of this iterative method: \Cref{eqn:predictor} makes a coarse prediction of the output, then~\cref{eqn:corrector} improves this prediction with a correction term derived from expanding the original PDE. This choice of $ \varphi $ falls within an ensemble viewpoint that the RFM may be used to improve pre-existing surrogate models by taking $ \varphi(\cdot;\theta) $ to be an existing emulator, but randomized in a principled way through $ \theta\sim\mu $.

For this particular example, we are cognizant of the facts that the random feature map $ \varphi $ requires full knowledge of the Darcy equation and a na\"ive evaluation of $ \varphi $ may be as expensive as solving the original PDE, which is itself a linear PDE; however, we believe that the ideas underlying the random features used here are intuitive and suggestive of what is possible in other applications areas. For example, RFMs may be applied on larger domains with simple geometries, viewed as supersets of the physical domain of interest, enabling the use of efficient algorithms such as the fast Fourier transform (FFT) even though these may not be available on the original problem, either because the operator to be inverted is spatially inhomogeneous or because of the complicated geometry of the physical domain. 

\section{Numerical Experiments}
\label{sec:experiment}
We now assess the performance of our proposed methodology on the approximation of operators $ \Fd:\cX\to\cY $ presented in~\Cref{sec:application}.
Practical implementation of the approach on a computer necessitates discretization of the input-output function spaces $ \cX $, $ \cY $. Hence in the numerical experiments that follow, all infinite-dimensional objects such as the training data, evaluations of random feature maps, and random fields are discretized on an equispaced mesh with $ K $ grid points to take advantage of the $ O(K\log K) $ computational speed of the FFT. The simple choice of equispaced points does not limit the proposed approach, as our formulation of the RFM on function space allows the method to be implemented numerically with any choice of spatial discretization. Such a numerical discretization procedure leads to the problem of high- but finite-dimensional approximation of discretized target operators mapping $ \R^{K} $ to $ \R^{K} $ by similarly discretized RFMs. However, we emphasize the fact that $ K $ is allowed to vary, and we study the properties of the discretized RFM as $ K $ varies, noting that since the RFM is defined conceptually on function space in \Cref{sec:problem} without reference to discretization, its discretized numerical realization has approximation quality consistent with the infinite-dimensional limit $ K\to\infty $. This implies that the same trained model can be deployed across the entire hierarchy of finite-dimensional spaces $ \R^{K} $ parametrized by $ K\in\N $ without the need to be re-trained, provided
$K$ is sufficiently large. Thus in this section, our notation does not make explicit the dependence of the discretized RFM or target operators on mesh size $ K $. We demonstrate these claimed properties numerically.

The input functions and our chosen random feature maps~\cref{eqn:rf_fourier} and~\cref{eqn:rf_predictor_corrector} require i.i.d. draws of Gaussian random fields to be fully defined. We efficiently sample these fields by truncating a Karhunen-Lo\'eve expansion and employing fast summation of the eigenfunctions with FFT. More precisely, on a mesh of size $ K $, denote by $ g(\cdot) $ a numerical approximation of a Gaussian random field on domain $ D=(0,1)^{d} $, $ d=1,\, 2 $:
\begin{equation}\label{eqn:grf_klexpansion}
g=\sum_{k\in Z_K}\xi_k\sqrt{\lambda_k}\phi_k \approx \sum_{k'\in\Z_{\geq 0}^{d}}\xi_{k'}\sqrt{\lambda_{k'}}\phi_{k'}  \sim N(0, C)\, ,
\end{equation}
where $ \{\xi_{j}\}\sim N(0,1) $ i.i.d. and $ Z_K\subset\Z_{\geq 0} $ is a truncated one-dimensional lattice of cardinality $ K $ ordered such that $ \{\lambda_j\} $ is non-increasing. The pairs $ (\lambda_{k'}, \phi_{k'}) $ are found by solving the eigenvalue problem $ C\phi_{k'}=\lambda_{k'}\phi_{k'} $ for non-negative, symmetric, trace-class operator $ C $ \cref{eqn:prior_covariance}. Concretely, these solutions are given by
\begin{equation}\label{eqn:eig_neumann}
\phi_{k'}(x)
=
\begin{cases}
\sqrt{2}\cos(k_1'\pi x_1)\cos(k_2'\pi x_2),& k_1' \ \text{or} \ k_2'=0\\
2\cos(k_1'\pi x_1)\cos(k_2'\pi x_2),& \text{otherwise} 
\end{cases}\, ,\quad
\lambda_{k'}=\tau^{2\al - 2}(\pi^{2}\abs{k'}^{2}+\tau^{2})^{-\al}
\end{equation}
for homogeneous Neumann boundary conditions when $ d=2 $, $ k'=(k_1', k_2')\in \Z_{\geq 0}^{2}{\setminus}\{0\} $, $ x=(x_1,x_2)\in (0,1)^{2} $, and given by
\begin{subequations}\label{eqn:periodic}
	\begin{align}
	\phi_{2j}(x)&=\sqrt{2}\cos(2\pi j x)\, ,\quad  \phi_{2j-1}(x)=\sqrt{2}\sin(2\pi j x)\, ,\quad \phi_0(x)=1\,,\\
	\lambda_{2j}&=\lambda_{2j-1}=\tau^{2\al - 1}(4\pi^{2}j^{2}+\tau^{2})^{-\al}\, ,\quad \lambda_0=\tau^{-1}
	\end{align}
\end{subequations}
for periodic boundary conditions when $ d=1 $, $ j\in\Z_{>0} $, and $ x\in (0,1) $. In both cases, we enforce that $ g $ integrate to zero over $ D $ by manually setting to zero the Fourier coefficient corresponding to multi-index $ k'=0 $. We use such $ g $ in all experiments that follow. Additionally, the $ k $ and $ k' $ used in this section to denote wavenumber indices should not be confused with our previous notation for kernels.

With the discretization and data generation setup now well-defined, and the pairs $ (\varphi, \mu) $ given in \Cref{sec:application}, the last algorithmic step is to train the RFM by solving~\cref{eqn:opt_normaleqn} and then test its performance. For a fixed number of random features $ m $, we only train and test a single realization of the RFM, viewed as a random variable itself. In each instance $ m $ is varied in the experiments that follow, the draws $ \{\theta_j\}_{j=1}^m $ are re-sampled i.i.d. from $ \mu $. 
To measure the distance between the trained RFM $ F_{m}(\cdot;\hat{\al}) $ and the ground truth map $ \Fd $, we employ the \emph{approximate expected relative test error}
\begin{equation}\label{eqn:test_error}
e_{n',m}\defby \dfrac{1}{n'}\sum_{j=1}^{n'}\dfrac{\norm{\Fd(a_j')-F_{m}(a_j';\hat{\al})}_{L^2}}{\norm{\Fd(a_j')}_{L^2}} \approx \E^{a'\sim\nu}\left[ \dfrac{\norm{\Fd(a')-F_{m}(a';\hat{\al})}_{L^2}}{\norm{\Fd(a')}_{L^2}}\right]\, ,
\end{equation}
where the $\{a_j'\}_{j=1}^{n'}$ are drawn i.i.d. from $\nu$ and $ n' $ denotes the number of input-output pairs used for testing. All $ L^2(D;\R) $ norms on the physical domain are numerically approximated by composite trapezoid rule quadrature. Since $ \cY\subset L^{2} $ for both the PDE solution operators \cref{eqn:solnmap_burg} and \cref{eqn:solnmap_darcy}, we also perform all required inner products during training in $ L^{2} $ rather than in $ \cY $; this results in smaller relative test error $ e_{n',m} $.

\subsection{Burgers' Equation: Experiment}\label{sec:burg_exp}
We generate a high resolution dataset of input-output pairs by solving Burgers' equation~\cref{eqn:burgers_ibvp} on an equispaced periodic mesh of size $ K=1025 $ (identifying the first mesh point with the last) with random initial conditions sampled from $ \nu=N(0,C) $ using \cref{eqn:grf_klexpansion}, where $ C $ is given 
by~\cref{eqn:prior_covariance} with parameter choices $\tau=7$ and $\al=2.5$. The full-order solver is a FFT-based pseudospectral method for spatial discretization \cite{fornberg1998practical} and a fourth-order Runge-Kutta integrating factor time-stepping scheme for time discretization \cite{kassam2005fourth}. All data represented on mesh sizes $ K<1025 $ used in both training and testing phases are subsampled from this original dataset, and hence we consider numerical realizations of $ \Fd $ \cref{eqn:solnmap_burg} up to $ \R^{1025}\to\R^{1025} $. We fix $ n=512 $ training and $ n'=4000 $ testing pairs unless otherwise noted, and also fix the viscosity to $ \ep=10^{-2} $ in all experiments. Lowering $ \ep $ leads to smaller length scale solutions and more difficult reconstruction; more data (higher $ n $) and features (higher $ m $) or a more expressive choice of $ (\varphi,\mu) $ would be required to achieve comparable error levels due to the slow decaying Kolmogorov width of the solution map. For simplicity, we set the forcing $ f\equiv 0 $, although nonzero forcing could lead to other interesting solution maps such as $ f\mapsto u(T,\cdot) $. It is easy to check that the solution will have zero mean for all time and a steady state of zero. Hence, we choose $ T\leq 2 $ to ensure that the solution is far enough away from steady state. For the random feature map~\cref{eqn:rf_fourier}, we fix the hyperparameters $ \al'=2 $, $ \tau'=5 $, $ \delta=0.0025 $, and $ \beta=4 $. The map itself is evaluated efficiently with the FFT and requires no other tools to be discretized. RFM hyperparameters were hand-tuned but not optimized. We find that regularization during training had a negligible effect for this problem, so the RFM is trained with $ \lambda=0 $ by solving the normal equations~\cref{eqn:opt_normaleqn} with the pseudoinverse to deliver the minimum norm least squares solution; we use the truncated SVD implementation in Python's $ \texttt{scipy.linalg.pinv2} $ for this purpose.
\begin{figure}[!htbp]
	\centering
	\begin{subfigure}[]{0.49\textwidth}
		\centering
		\includegraphics[width=\textwidth]{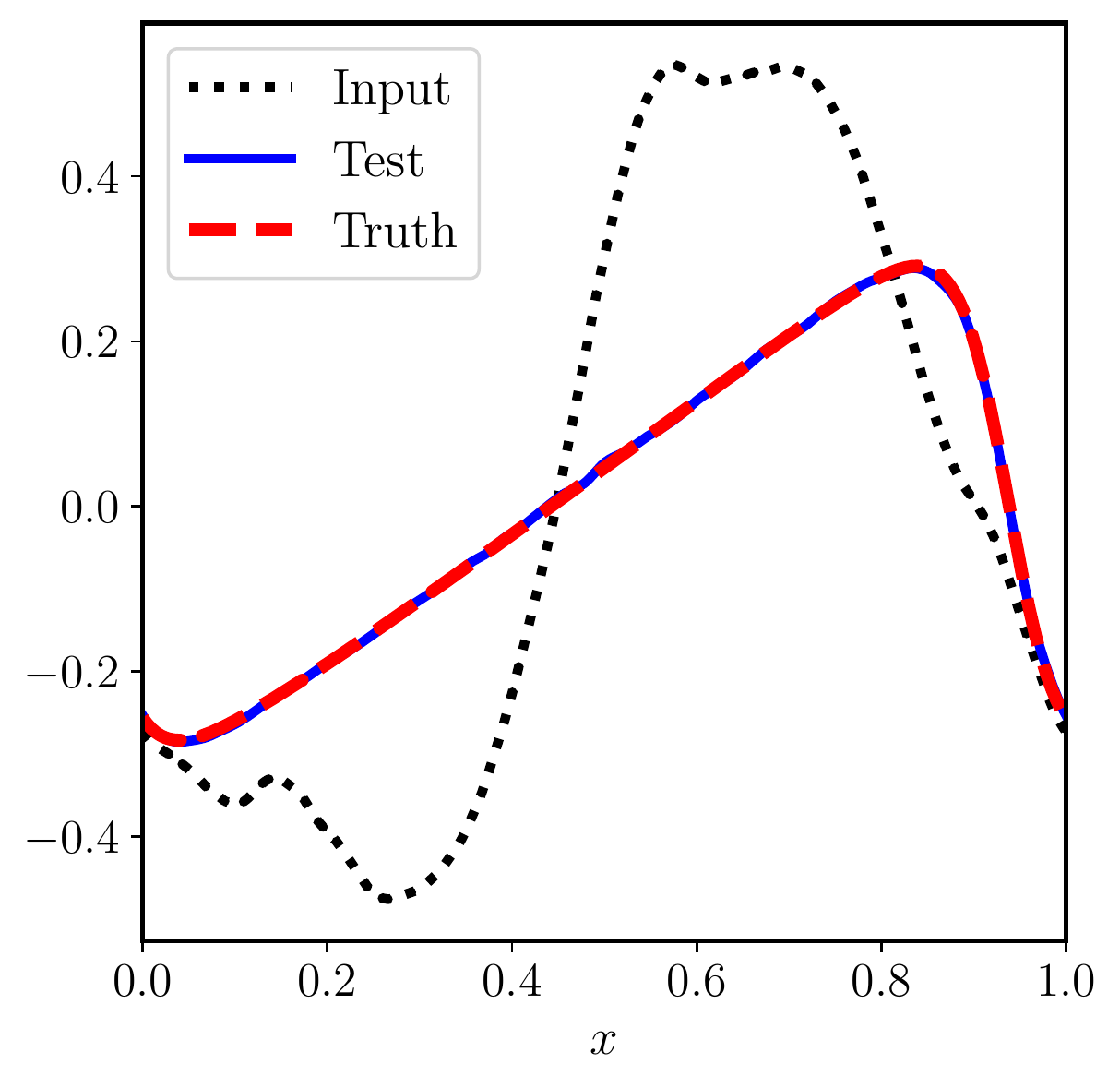}
		\caption{}
		\label{fig:prediction_onesample}
	\end{subfigure}%
	\hfill%
	\begin{subfigure}[]{0.49\textwidth}
		\centering
		\includegraphics[width=\textwidth]{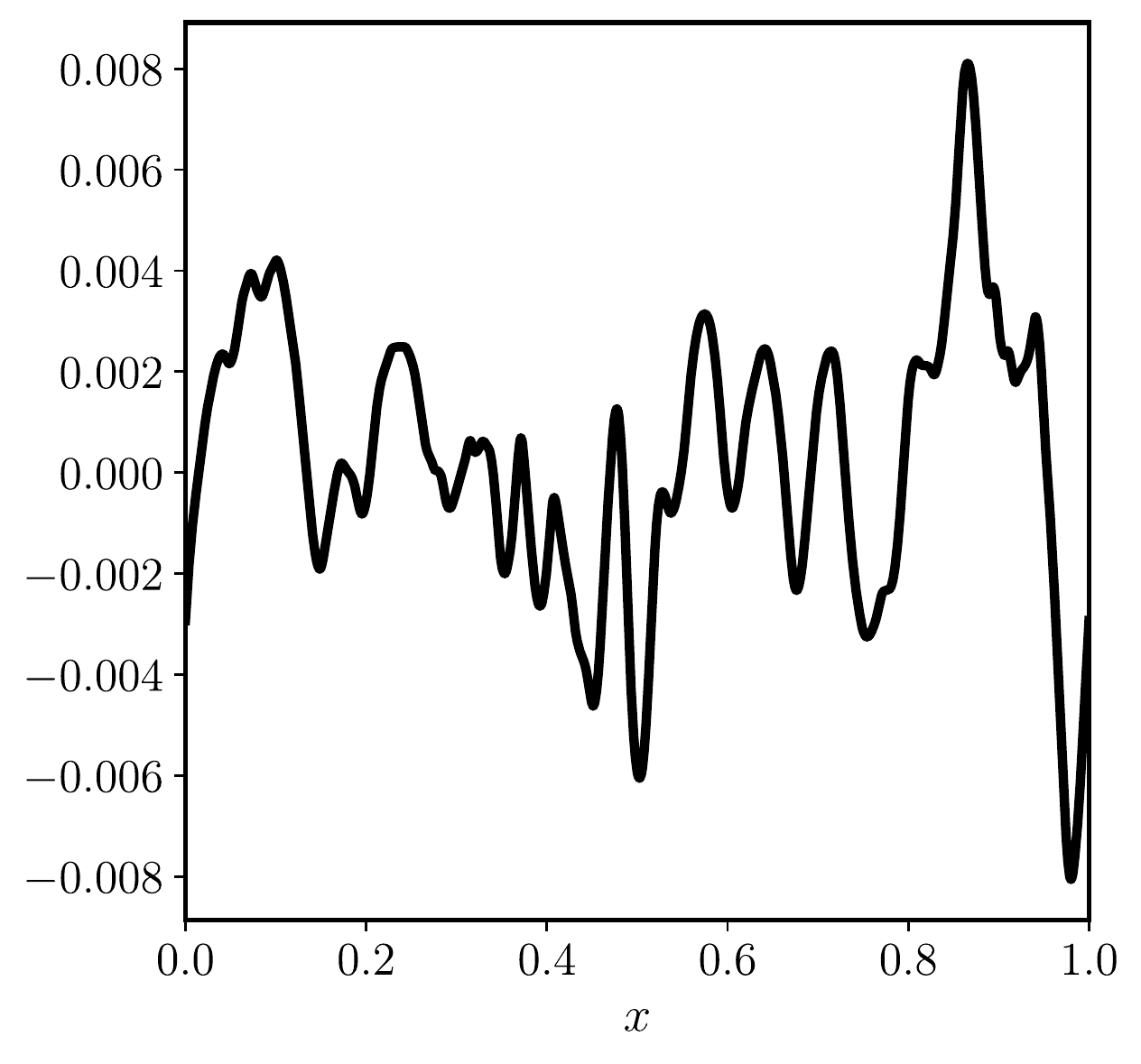}
		\caption{}
		\label{fig:pwerror_onesample}
	\end{subfigure}
	\vspace{-5mm}
	\caption{Representative input-output test sample for the Burgers' equation solution map $ \Fd\defby \Psi_{1} $: Here, $ n=512 $, $ m=1024 $, and $ K=1025 $. \Cref{fig:prediction_onesample} shows a sample input, output (truth), and trained RFM prediction (test), while \cref{fig:pwerror_onesample} displays the pointwise error. The relative $ L^2 $ error for this single prediction is $ 0.0146 $.}
	\label{fig:sample_burg}
\end{figure}

Our experiments study the RFM approximation to the viscous Burgers' equation evolution operator semigroup~\cref{eqn:solnmap_burg}. As a visual aid for the high-dimensional problem at hand, \Cref{fig:sample_burg} shows a representative sample input and output along with a trained RFM test prediction. To determine whether the RFM has actually learned the correct evolution operator, we test the semigroup property of the map;~\cite{wu2020data} pursues closely related work also in a Fourier space setting. Denote the $ (j-1) $-fold composition of a function $ G $ with itself by $ G^{j} $. Then, with $ u(0,\cdot)=a $, we have
\begin{equation}\label{eqn:semigroup}
(\Psi_{T}\circ \cdots \circ \Psi_{T})(a)=\Psi_{T}^{j}(a)=\Psi_{jT}(a)=u(jT, \cdot)
\end{equation}
by definition. We train the RFM on input-output pairs from the map $ \Psi_{T} $ with $ T\defby 0.5 $ to obtain $ \hat{F}\defby F_{m}(\cdot;\hat{\al}) $. Then, it should follow from~\cref{eqn:semigroup} that $ \hat{F}^{j}\approx \Psi_{jT}$, that is, each application of $ \hat{F} $ should evolve the solution $ T $ time units. We test this semigroup approximation by learning the map $ \hat{F} $ and then comparing $ \hat{F}^{j} $ on  $n'= 4000$ fixed inputs to outputs from each of the operators $ \Psi_{jT} $, with $ j\in\{1,2,3,4\} $ (the solutions at time $ T $, $ 2T $, $ 3T $, $ 4T $). The results are presented in~\Cref{tab:time_upscale} for a fixed mesh size $ K=129 $.
\begin{table}[!htbp]
	\centering
	\begin{tabular}{ @{}lccccc@{} }
		\toprule
		Train on: & $ T=0.5 $ & $\ \quad \ $ Test on: & $ 2T=1.0 $ & $ 3T=1.5 $ & $4T=2.0$ \\
		\midrule
		& 0.0360 & & 0.0407 & 0.0528 & 0.0788 \\
		\bottomrule
	\end{tabular}
	\vspace{1mm}
	\caption{Expected relative error $ e_{n',m} $ for time upscaling with the learned RFM operator semigroup for Burgers' equation: Here, $ n'=4000 $, $ m=1024 $, $ n=512 $, and $ K=129 $. The RFM is trained on data from the evolution operator $ \Psi_{T=0.5} $, and then tested on input-output samples generated from $ \Psi_{jT} $, where $ j=2,\, 3,\, 4 $, by repeated composition of the learned model. The increase in error is small even after three compositions, reflecting excellent out-of-distribution performance.}
	\label{tab:time_upscale}
\end{table}
We observe that the composed RFM map $ \hat{F}^{j} $ accurately captures $ \Psi_{jT} $, though this accuracy deteriorates as $ j $ increases due to error propagation in time as is common with any traditional integrator. However, even after three compositions corresponding to 1.5 time units past the training time $ T=0.5 $, the relative error only increases by around $ 0.04 $. It is remarkable that the RFM learns time evolution without explicitly time-stepping the PDE~\cref{eqn:burgers_ibvp} itself. Such a procedure is coined \emph{time upscaling} in the PDE context and in some sense breaks the CFL stability barrier~\cite{demanet2006curvelets}. \Cref{tab:time_upscale} is evidence that the RFM has excellent out-of-distribution performance: although only trained on inputs $ a\sim\nu $, the model outputs accurate predictions given new input 
samples $ \Psi_{jT}(a) \sim (\Psi_{jT})_{\sharp}\nu $.

\begin{figure}[!htbp]
	\centering
	\begin{subfigure}[]{0.49\textwidth}
		\centering
		\includegraphics[width=\textwidth]{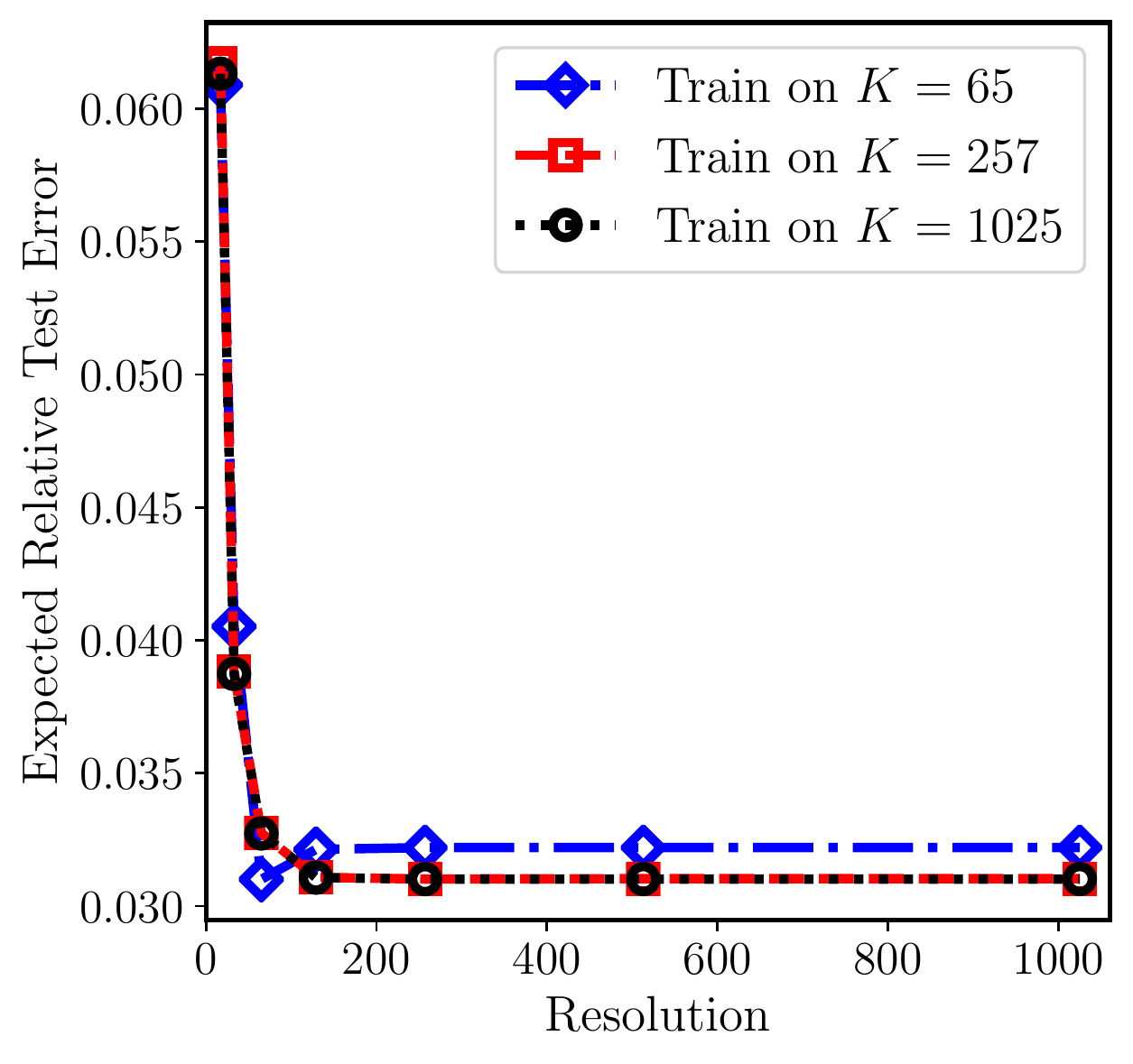}
		\caption{}
		\label{fig:gridtransfer_burg_panel}
	\end{subfigure}%
	\hfill%
	\begin{subfigure}[]{0.49\textwidth}
		\centering
		\includegraphics[width=\textwidth]{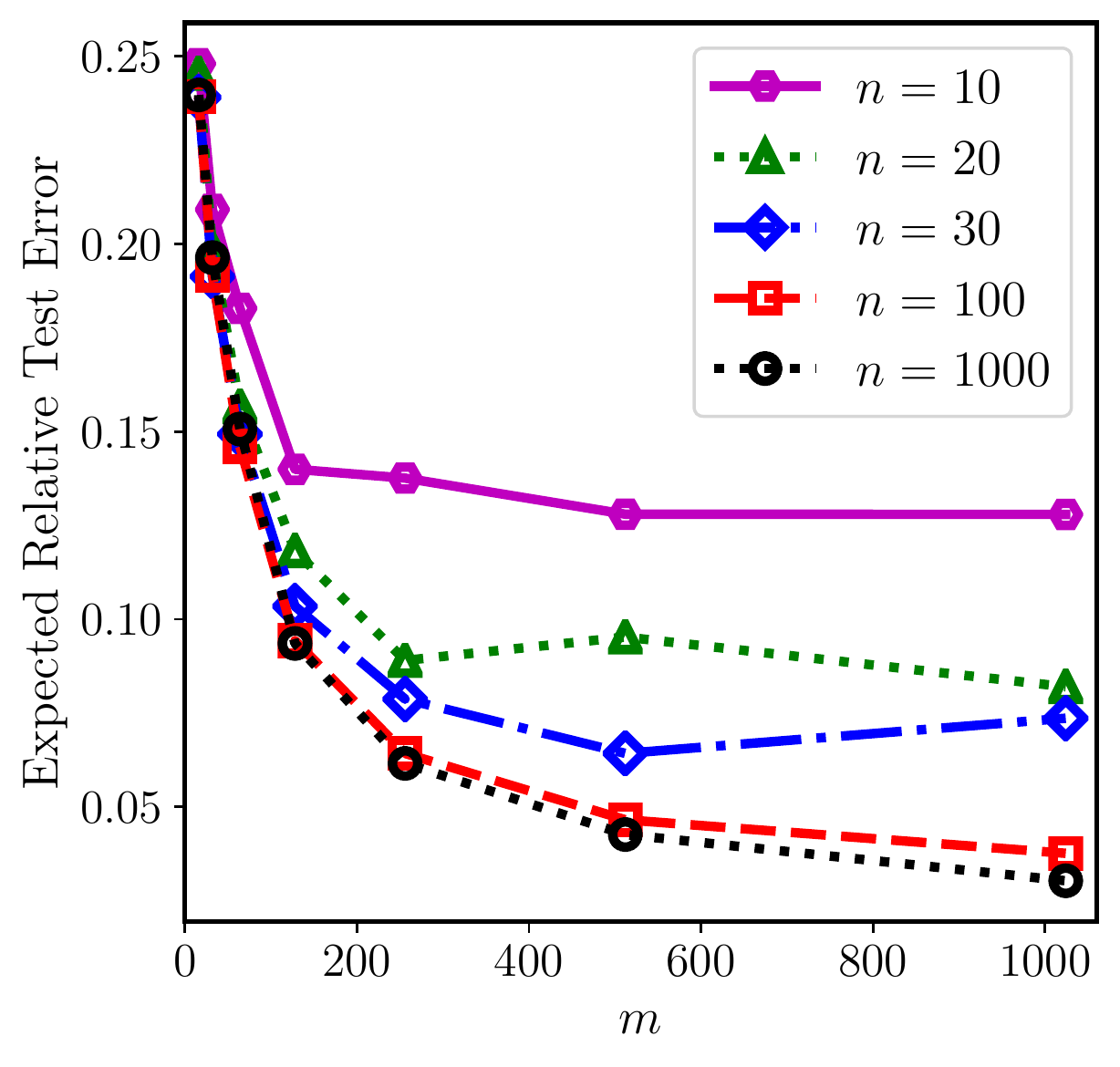}
		\caption{}
		\label{fig:gridsweep_burg_n}
	\end{subfigure}
	\vspace{-5mm}
	\caption{Expected relative test error of a trained RFM for the Burgers' evolution operator $ \Fd=\Psi_{1} $ with $ n'=4000 $ test pairs: \Cref{fig:gridtransfer_burg_panel} displays the invariance of test error w.r.t. training and testing on different resolutions for $ m=1024 $ and $ n=512 $ fixed; the RFM can train and test on different mesh sizes without loss of accuracy. \Cref{fig:gridsweep_burg_n} shows the decay of the test error for resolution $ K=129 $ fixed as a function of $ m $ and $ n $; the smallest error achieved is $ 0.0303 $ for $ n=1000 $ and $ m=1024 $.}
	\label{fig:gridtranfser_burg}
\end{figure}

We next study the ability of the RFM to transfer its learned coefficients $ \hat{\al} $ obtained from training on mesh size $ K $ to different mesh resolutions $ K' $ in~\cref{fig:gridtransfer_burg_panel}. We fix $ T\defby 1 $ from here on and observe that the lowest test error occurs when $ K=K'$, that is, when the train and test resolutions are identical; this behavior was also observed in the contemporaneous work~\cite{li2020neural}. At very low resolutions, such as $ K=17 $ here, the test error is dominated by discretization error which can become quite large; for example, resolving conceptually infinite-dimensional objects such as the Fourier space-based feature map in \cref{eqn:rf_fourier} or the $ L^2 $ norms in \cref{eqn:test_error} with only $ 17 $ grid points gives bad accuracy. But outside this regime, the errors are essentially constant across resolution regardless of the training resolution $ K $, indicating that the RFM learns its optimal coefficients independently of the resolution and hence generalizes well to any desired mesh size. In fact,
the trained model could be deployed on different discretizations of the domain $ D $ (e.g.: various choices of finite elements, graph-based/particle methods), not just with different mesh sizes. Practically speaking, this means that high resolution training sets can be subsampled to smaller mesh sizes $ K $ (yet still large enough to avoid large discretization error) for faster training, leading to a trained model with nearly the same accuracy at all higher resolutions.

The smallest expected relative test error achieved by the RFM is $ 0.0303 $ for the configuration detailed in \cref{fig:gridsweep_burg_n}. This excellent performance is encouraging because the error we report is of the same order of magnitude as that reported in Sec. 5.1 of \cite{li2020fourier} for the same Burgers' solution operator that we study, but with slightly different problem parameter choices. We emphasize that the Neural Operator methods in that work are based on deep learning, which involves training neural networks by solving a non-convex optimization problem with stochastic gradient descent, while our random feature methods have orders of magnitude fewer trainable parameters that are easily optimized through convex optimization. In \cref{fig:gridsweep_burg_n}, we also note that for a small number of training data $ n $, the error does not always decrease as the number of random features $ m $ increases. This indicates a delicate dependence of $ m $ as a function of $ n $, in particular, $ n $ must increase with $ m $ as is expected from parametric estimation; we observe the desired monotonic decrease in error with $ m $ when $ n $ is increased to $ 100 $ or $ 1000 $. In the over-parametrized regime, the authors in~\cite{ma2019generalization} present a loose bound for this dependence for real-valued outputs. We leave a detailed account of the dependence of $ m $ on $ n $ required to achieve a certain error tolerance to future work and refer the interested reader to~\cite{caponnetto2007optimal} for detailed statistical analysis in a related setting.

\begin{figure}[!htbp]
	\centering
	\begin{subfigure}[]{0.49\textwidth}
		\centering
		\includegraphics[width=\textwidth]{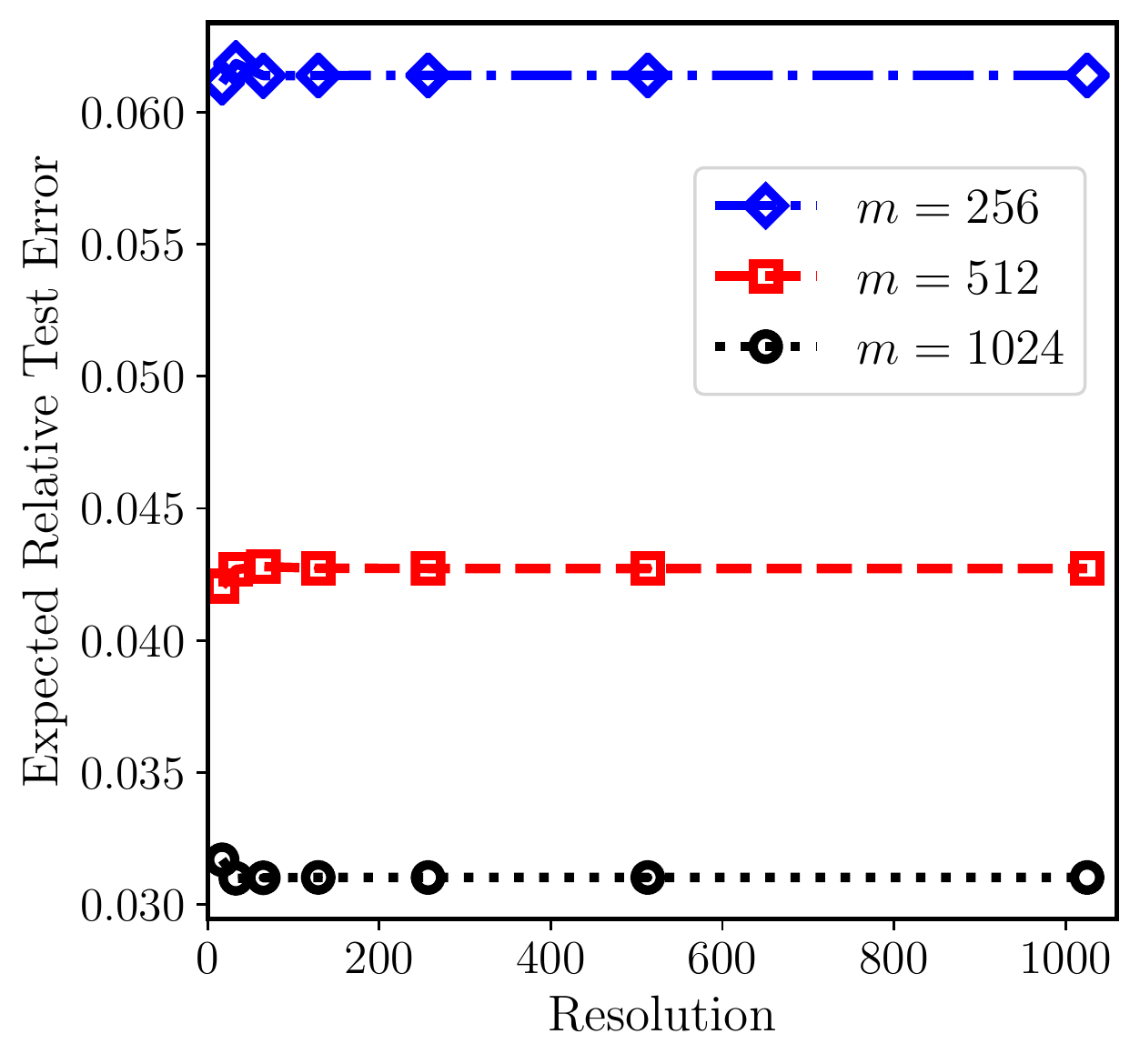}
		\caption{}
		\label{fig:gridsweep_burg1}
	\end{subfigure}%
	\hfill%
	\begin{subfigure}[]{0.49\textwidth}
		\centering
		\includegraphics[width=\textwidth]{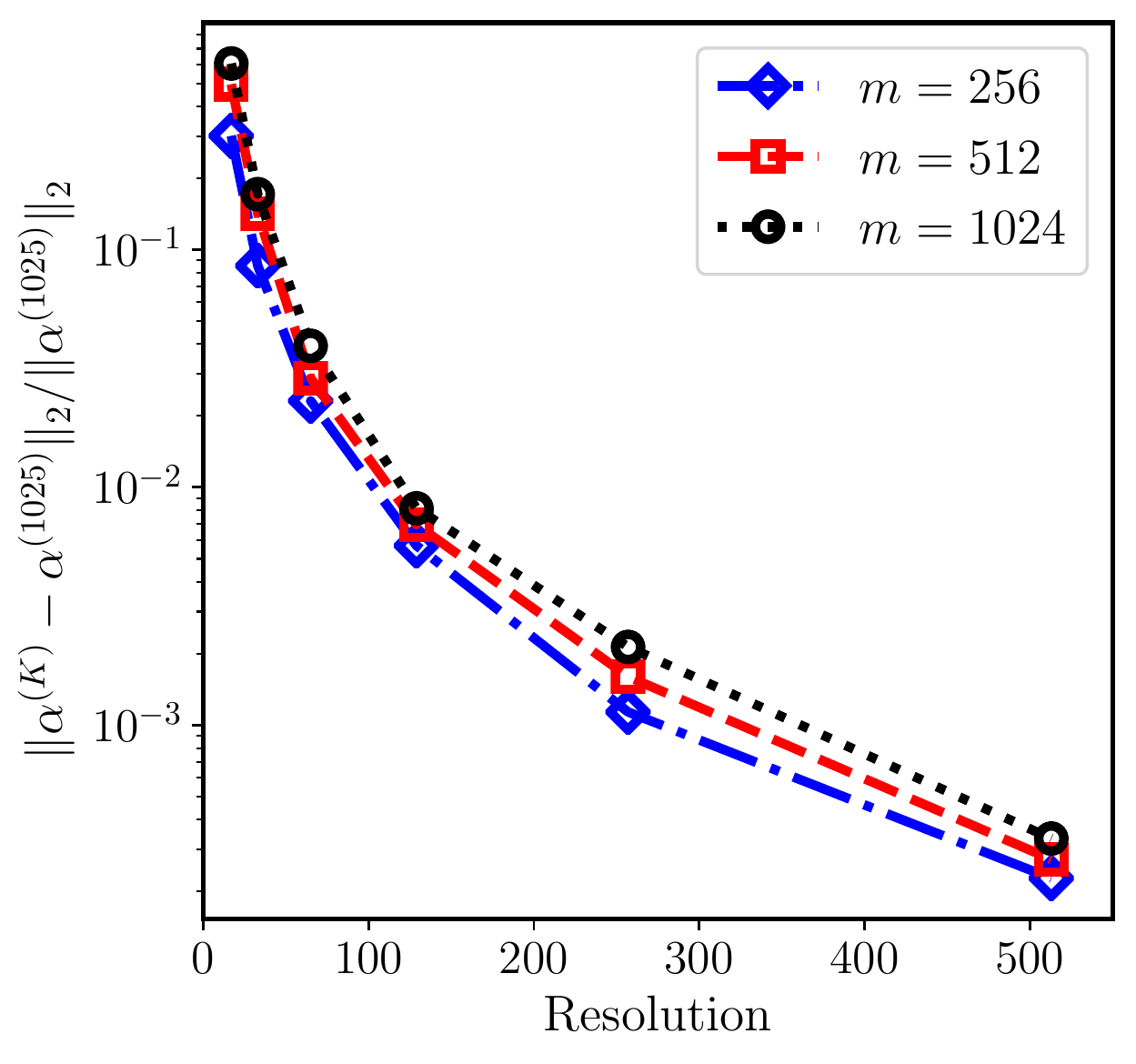}
		\caption{}
		\label{fig:gridsweep_burg2}
	\end{subfigure}
   	\vspace{-5mm}
	\caption{Results of a trained RFM for the Burgers' equation evolution operator $ \Fd=\Psi_{1} $: Here, $ n=512 $ training and $ n'=4000 $ testing pairs were used. \Cref{fig:gridsweep_burg1} shows resolution-invariant test error for various $ m $; the error follows the $ O(m^{-1/2}) $ Monte Carlo rate remarkably well. \Cref{fig:gridsweep_burg2} displays the relative error of the learned coefficient $ \al $ w.r.t. the coefficient learned on the highest mesh size ($ K=1025 $).}
	\label{fig:gridsweep_burg}
\end{figure}

Finally, \Cref{fig:gridsweep_burg} demonstrates the invariance of the expected relative test error to the mesh resolution used for training and testing. This result is a consequence of framing the RFM on function space; other machine learning-based surrogate methods defined in finite-dimensions exhibit an \emph{increase} in test error as mesh resolution is increased~(see \cite{bhattacharya2020pca}, Sec. 4, for a numerical account of this phenomenon). The first panel,~\Cref{fig:gridsweep_burg1}, shows the error as a function of mesh resolution for three values of $ m $. For very low resolution, the error varies slightly but then flattens out to a constant value as $ K\to\infty $. More interestingly, these constant values of error, $ e_{n',m}=0.063 $, $ 0.043 $, and $ 0.031 $ corresponding to $ m=256$, $512$, and $ 1024 $, respectively, closely match the Monte Carlo rate $ O(m^{-1/2}) $. While more theory is required to understand this behavior, it suggests that the optimization process finds coefficients close to those arising from a Monte Carlo approximation of $ \Fd $ as discussed in~\Cref{sec:rfm}. The second panel,~\Cref{fig:gridsweep_burg2}, indicates that the learned coefficient $ \al^{(K)} $ for each $ K $ converges to some $ \al^{(\infty)} $ as $ K\to\infty $, again reflecting the design of the RFM as a mapping between infinite-dimensional spaces.

\subsection{Darcy Flow: Experiment}\label{sec:darcy_exp}
In this section, we consider Darcy flow on the physical domain $ D\defby (0,1)^{2} $, the unit square. We generate a high resolution dataset of input-output pairs for $ \Fd $~\cref{eqn:solnmap_darcy} by solving~\Cref{eqn:darcy_flow} on an equispaced $ 257\times 257 $ mesh (size $ K=257^{2} $) using a second order finite difference scheme. All mesh sizes $ K<257^{2} $ are subsampled from this original dataset and hence we consider numerical realizations of $ \Fd $ up to $ \R^{66049}\to\R^{66049} $. We denote \emph{resolution} by $ r $ such that $ K=r^{2} $. We fix $ n=128 $ training and $ n'=1000 $ testing pairs unless otherwise noted. The input data are drawn from the level set measure $ \nu $~\cref{eqn:prior_levelset} with $ \tau=3 $ and $ \al=2 $ fixed. We choose $ a^{+}=12 $ and $ a^{-}=3 $ in all experiments that follow and hence the contrast ratio $ a^{+}/a^{-}=4 $ is fixed. The source is fixed to $ f\equiv 1 $, the constant function. We evaluate the predictor-corrector random features $ \varphi $~\cref{eqn:rf_predictor_corrector} using an FFT-based fast Poisson solver corresponding to an underlying second order finite difference stencil at a cost of $ O(K\log K) $ per solve. The smoothed coefficient $ a_{\ep} $ in the definition of $ \varphi $ is obtained by solving~\cref{eqn:heat_mollify_neumann} with time step $ 0.03 $ and diffusion constant $ \eta=10^{-4} $; with centered second order finite differences, this incurs 34 time steps and hence a cost $O(34K)$. We fix the hyperparameters $ \al'=2 $, $ \tau'=7.5 $, $ s^{+}=1/12 $, $ s^{-}=-1/3 $, and $ \delta=0.15 $ for the map $ \varphi $. Unlike in~\Cref{sec:burg_exp}, we find via grid search on $ \lambda $ that regularization during training does improve the reconstruction of the Darcy flow solution operator and hence we train with $ \lambda\defby 10^{-8} $ fixed. We remark that, for simplicity, the above hyperparameters were not systematically and jointly optimized; as a consequence the RFM performance has the capacity to improve beyond the results in this section.

\begin{figure}[!htbp]
	\centering
	\begin{subfigure}[]{0.49\textwidth}
		\centering
		\includegraphics[width=\textwidth]{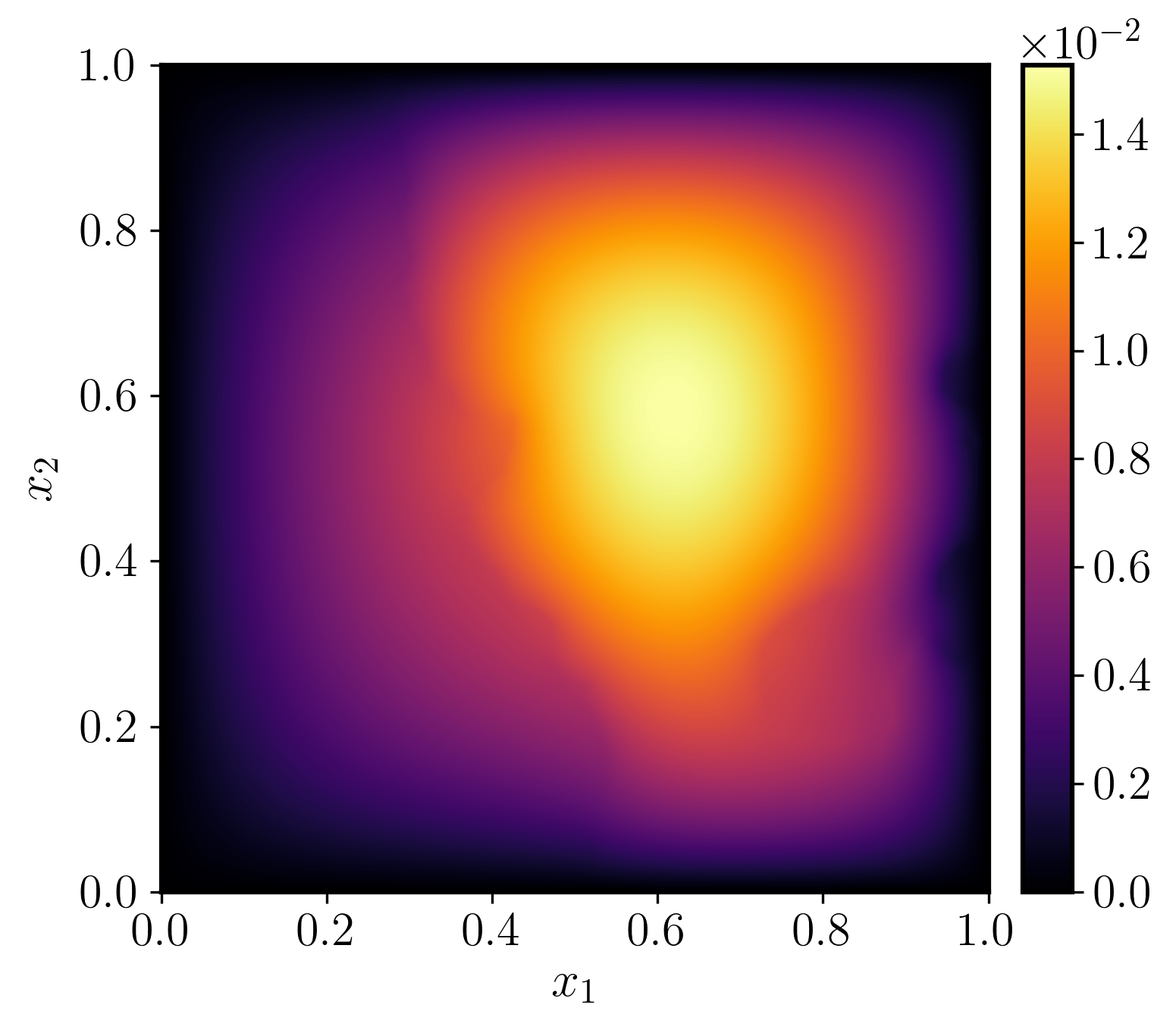}
		\caption{Truth}
		\label{fig:darcy_truth}
	\end{subfigure}%
	\hfill%
	\begin{subfigure}[]{0.49\textwidth}
		\centering
		\includegraphics[width=\textwidth]{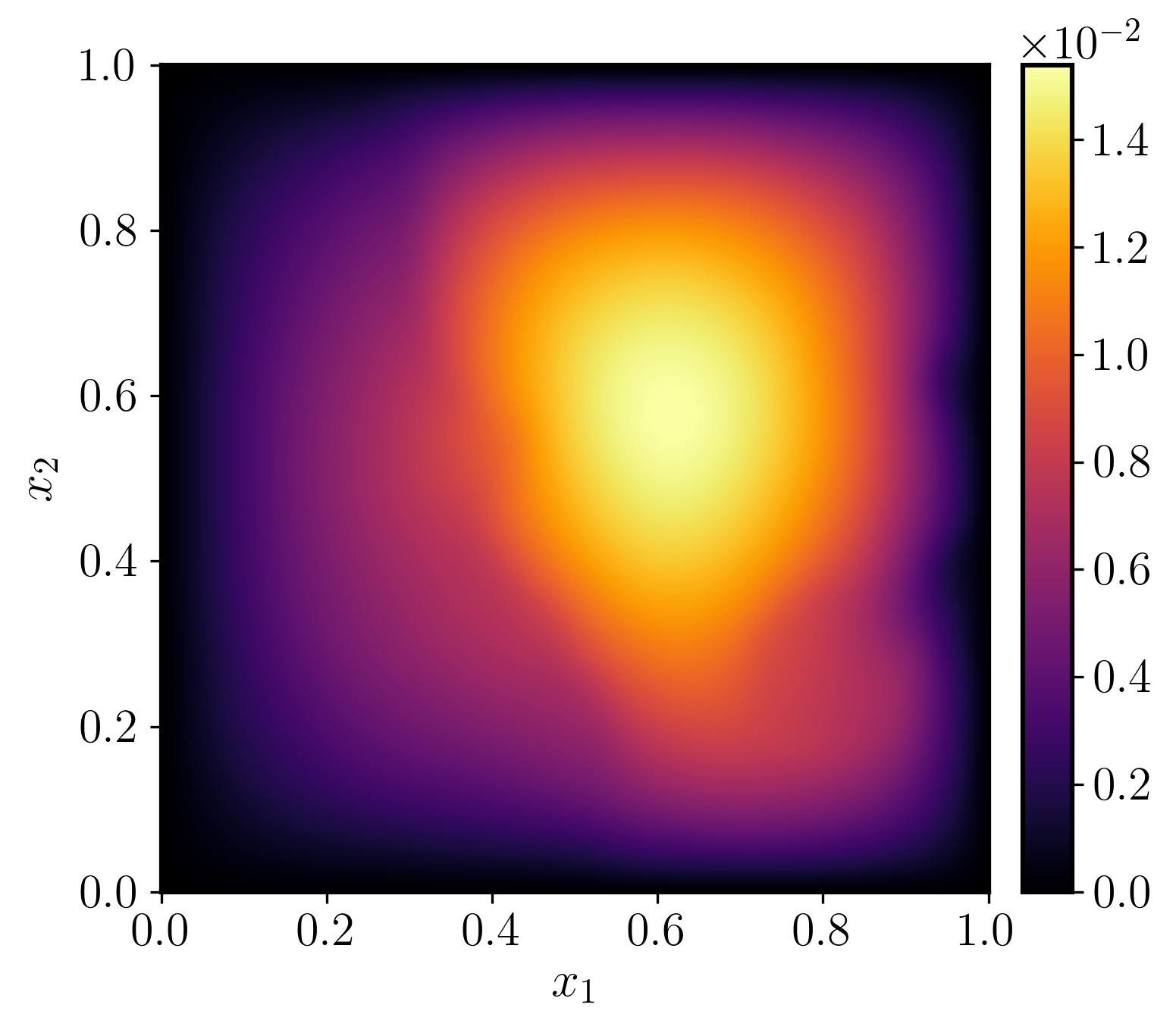}
		\caption{Approximation}
		\label{fig:darcy_predict}
	\end{subfigure}
	\begin{subfigure}[]{0.49\textwidth}
		\centering
		\includegraphics[width=\textwidth]{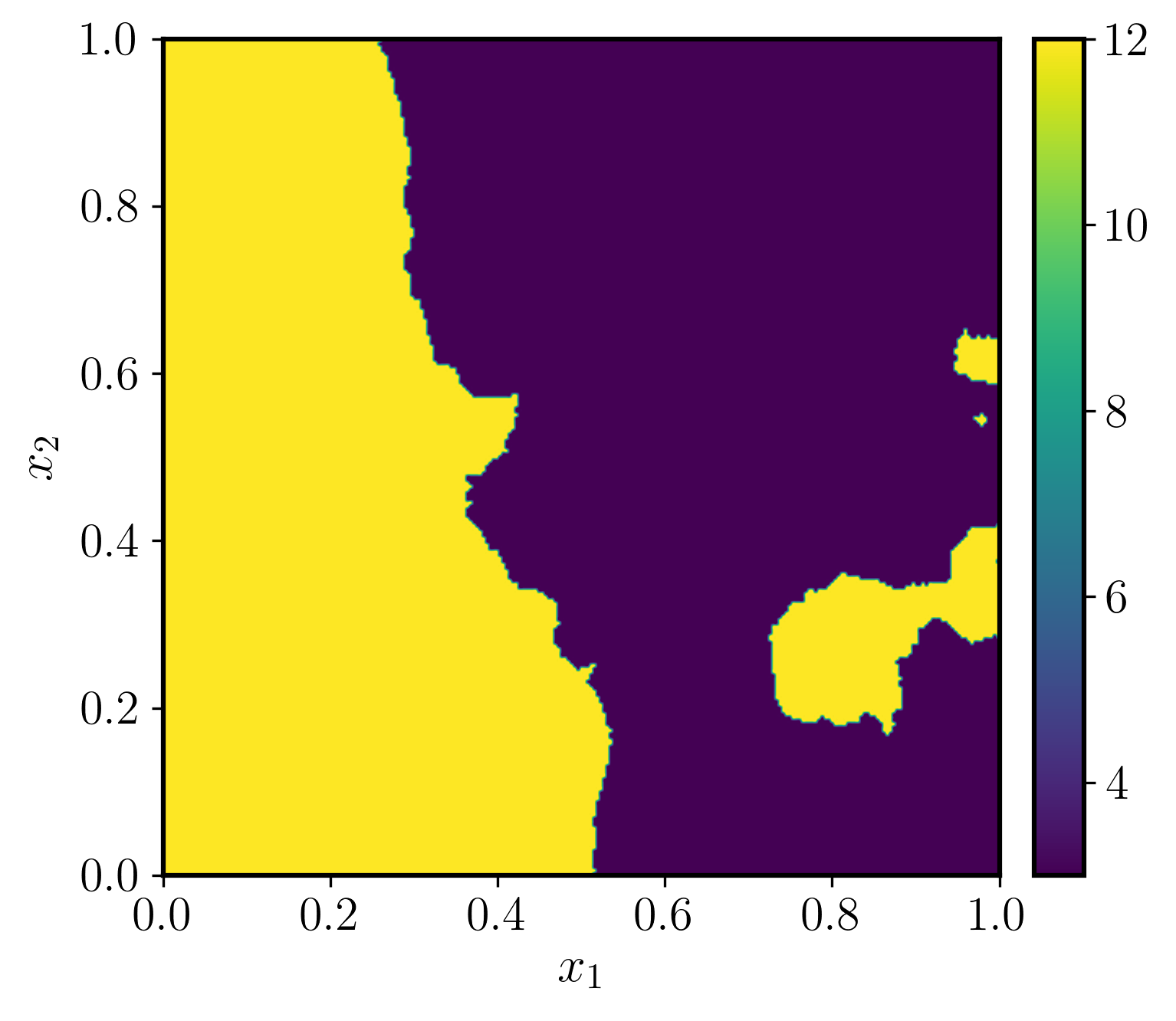}
		\caption{Input}
		\label{fig:darcy_input}
	\end{subfigure}%
	\hfill%
	\begin{subfigure}[]{0.49\textwidth}
		\centering
		\includegraphics[width=\textwidth]{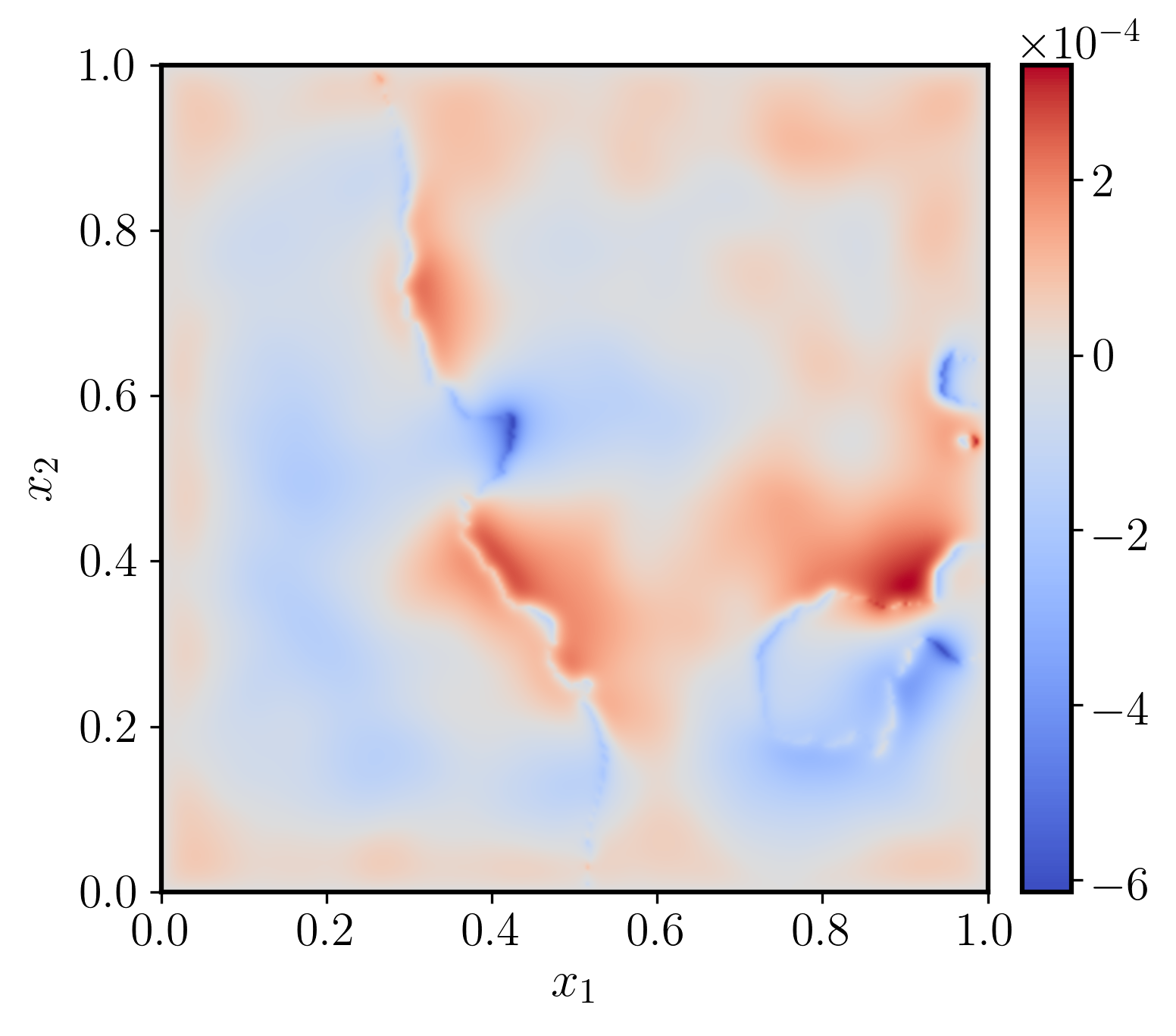}
		\caption{Pointwise Error}
		\label{fig:darcy_pwerror}
	\end{subfigure}
	\vspace{-5mm}
	\caption{Representative input-output test sample for the Darcy flow solution map: Here, $ n=256 $, $ m=350 $, and $ K=257^{2} $. \Cref{fig:darcy_input} shows a sample input, \Cref{fig:darcy_truth} the resulting output (truth), \Cref{fig:darcy_predict} a trained RFM prediction, and \Cref{fig:darcy_pwerror} the pointwise error. The relative $ L^2 $ error for this single prediction is $ 0.0122 $.}
	\label{fig:darcy_inputoutput}
\end{figure}

\begin{figure}[!htbp]
	\centering
	\begin{subfigure}[]{0.49\textwidth}
		\centering
		\includegraphics[width=\textwidth]{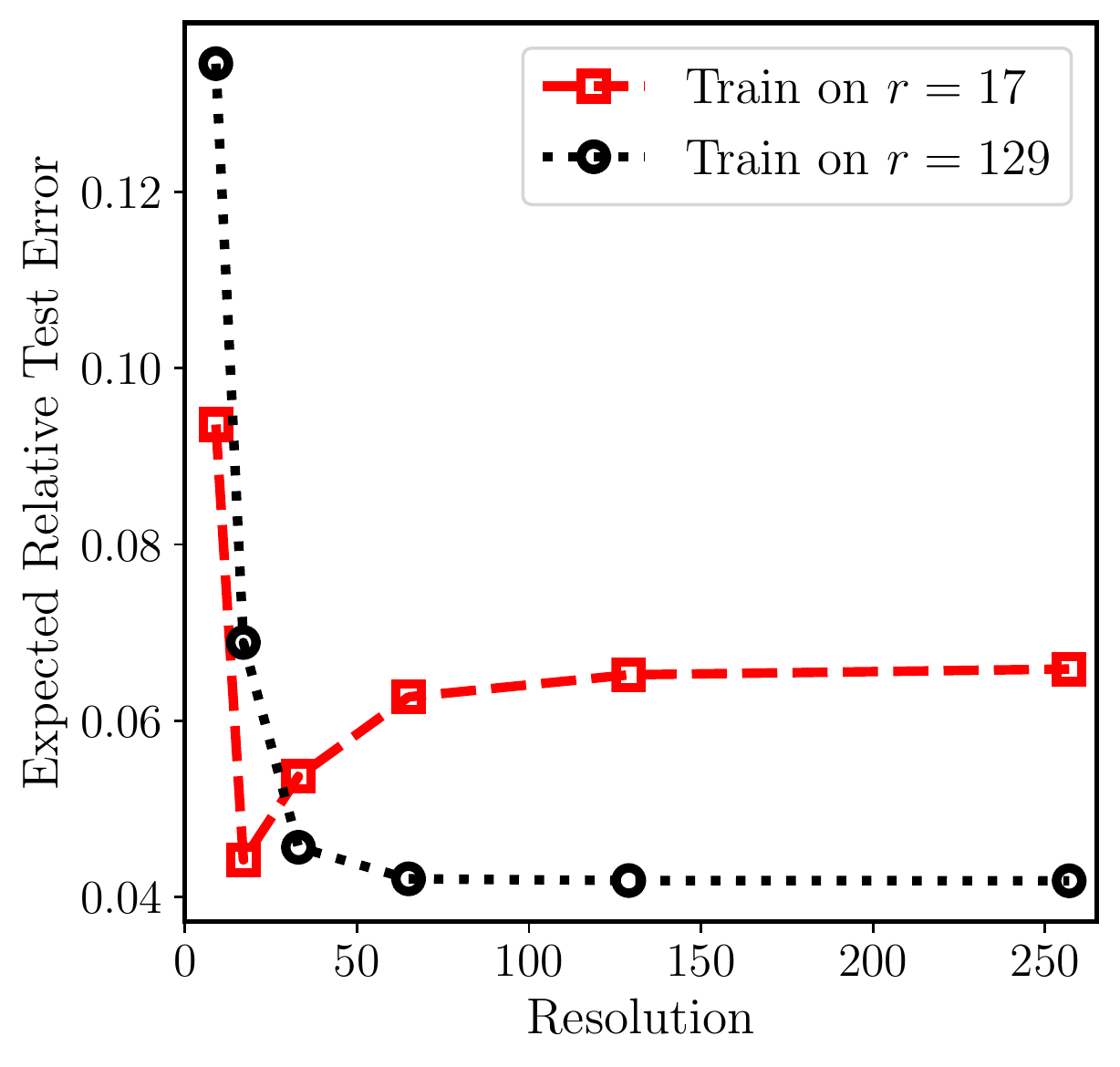}
		\caption{}
		\label{fig:gridtransfer_darcy_panel}
	\end{subfigure}%
	\hfill%
	\begin{subfigure}[]{0.49\textwidth}
		\centering
		\includegraphics[width=\textwidth]{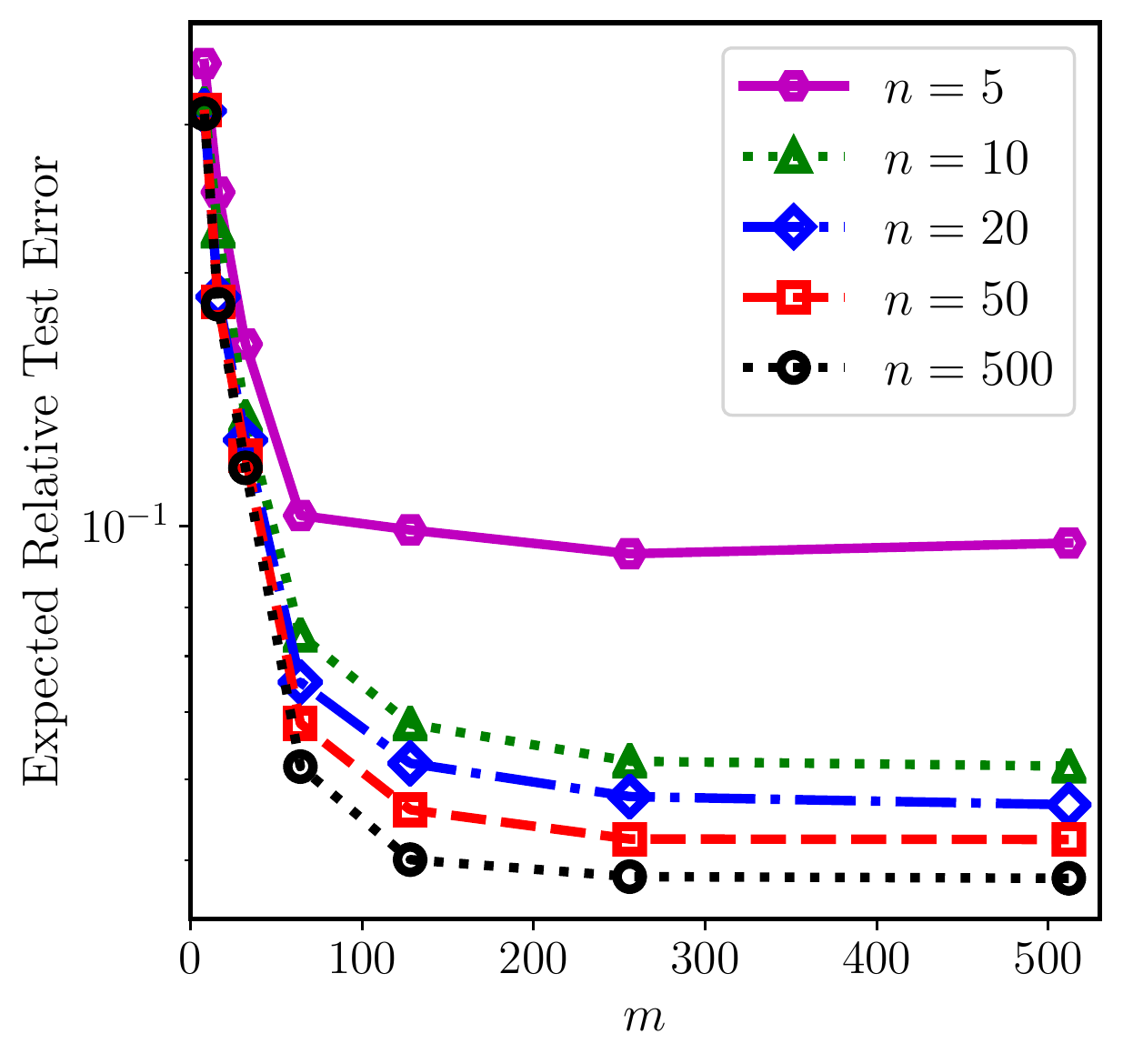}
		\caption{}
		\label{fig:gridsweep_darcy_n}
	\end{subfigure}
	\vspace{-5mm}
	\caption{Expected relative test error of a trained RFM for Darcy flow with $ n'= 1000 $ test pairs: \Cref{fig:gridtransfer_darcy_panel} displays the invariance of test error w.r.t. training and testing on different resolutions for $ m=512 $ and $ n=256 $ fixed; the RFM can train and test on different mesh sizes without significant loss of accuracy. \Cref{fig:gridsweep_darcy_n} shows the decay of the test error for resolution $ r=33 $ fixed as a function of $ m $ and $ n $; the smallest error achieved is $ 0.0381 $ for $ n=500 $ and $ m=512 $.}
	\label{fig:gridtranfser_darcy}
\end{figure}

Darcy flow is characterized by the geometry of the high contrast coefficients $ a\sim\nu $. As seen in \Cref{fig:darcy_inputoutput}, the solution inherits the steep interfaces of the input. However, we see that a trained RFM with predictor-corrector random features \cref{eqn:rf_predictor_corrector} captures these interfaces well, albeit with slight smoothing; the error concentrates on the location of the interface. The effect of increasing $ m $ and $ n $ on the test error is shown in~\Cref{fig:gridsweep_darcy_n}. Here, the error appears to saturate more than was observed for the Burgers' equation problem~(\Cref{fig:gridsweep_burg_n}). However, the smallest test error achieved for the best performing RFM configuration is $ 0.0381 $, which is on the same scale as the error reported in competing neural network-based methods~\cite{bhattacharya2020pca,li2020neural} for the same Darcy flow setup.

The RFM is able to be successfully trained and tested on different resolutions for Darcy flow. \Cref{fig:gridtransfer_darcy_panel} shows that, again, for low resolutions, the smallest relative test error is achieved when the train and test resolutions are identical (here, for $ r=17 $). However, when the resolution is increased away from this low resolution regime, the relative test error slightly increases then approaches a constant value, reflecting the function space design of the method. Training the RFM on a high resolution mesh poses no issues when transferring to lower or higher resolutions for model evaluation, and it achieves consistent error for test resolutions sufficiently large~(i.e., $ r \geq 33 $, the regime where discretization error starts to become negligible). Additionally, the RFM basis functions $ \{\varphi(\cdot;\theta_j)\}_{j=1}^{m} $ are defined without any dependence on the training data unlike in other competing approaches based on similar shallow linear approximations, such as the reduced basis method or the PCA-NN method in~\cite{bhattacharya2020pca}. Consequently, our random feature model may be directly evaluated on any desired mesh resolution once trained (``super-resolution''), whereas those aforementioned approaches require some form of interpolation to transfer between different mesh sizes (see \cite{bhattacharya2020pca}, Sec. 4.3).

\begin{figure}[!htbp]
	\centering
	\begin{subfigure}[]{0.49\textwidth}
		\centering
		\includegraphics[width=\textwidth]{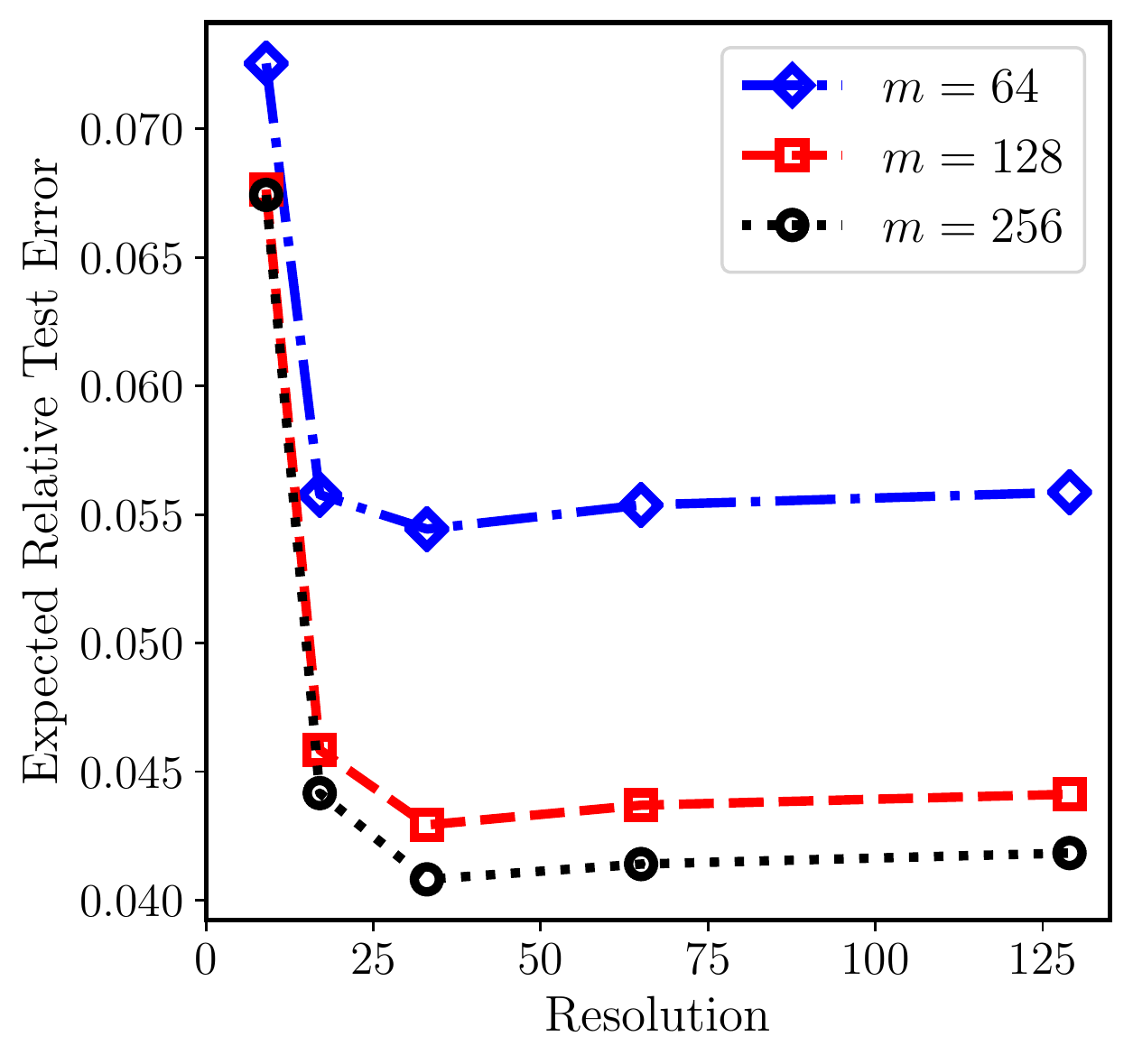}
		\caption{}
		\label{fig:gridsweep_darcy1}
	\end{subfigure}%
	\hfill%
	\begin{subfigure}[]{0.49\textwidth}
		\centering
		\includegraphics[width=\textwidth]{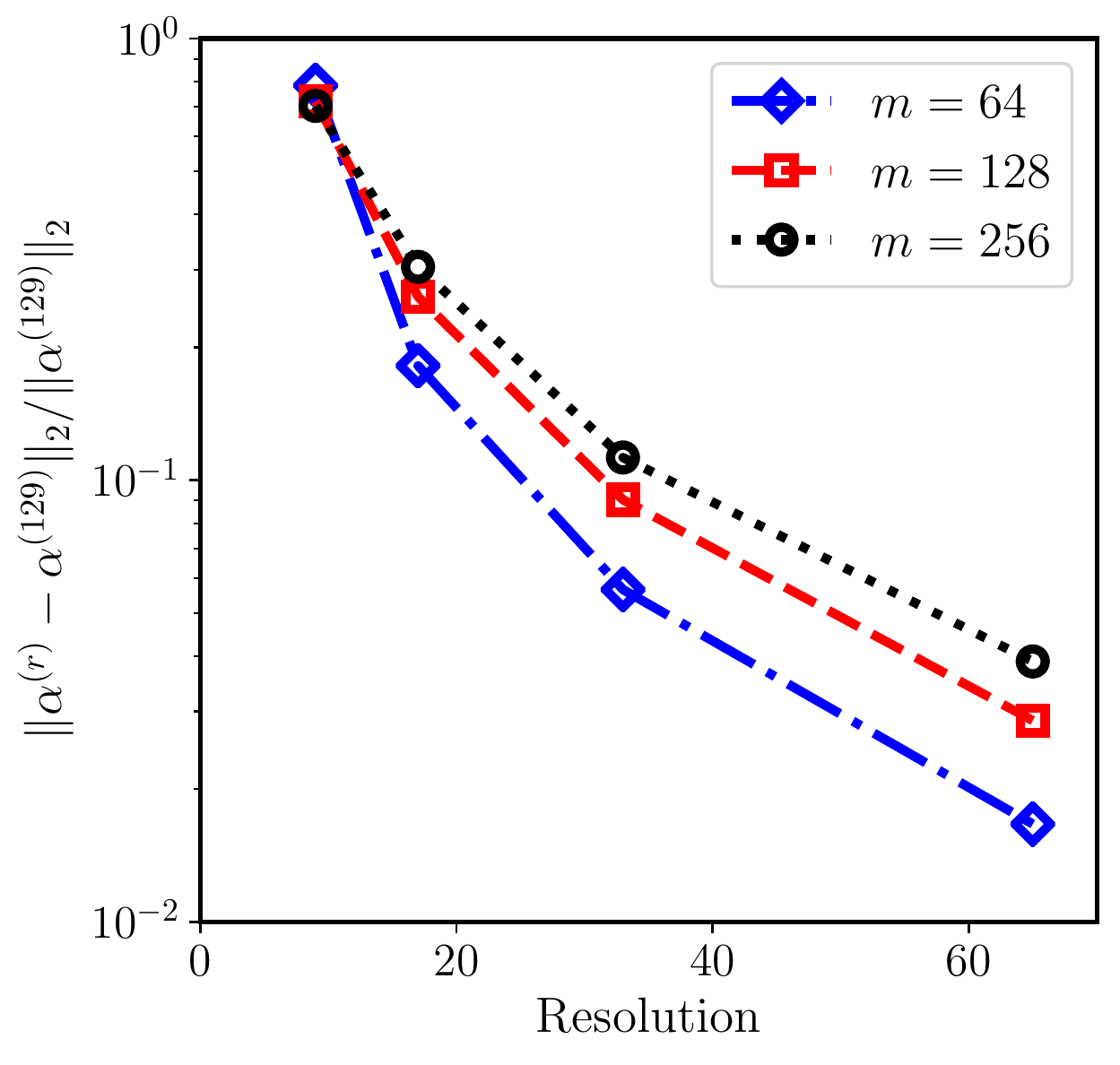}
		\caption{}
		\label{fig:gridsweep_darcy2}
	\end{subfigure}
	\vspace{-5mm}
	\caption{Results of a trained RFM for Darcy flow: Here, $ n=128 $ training and $ n'=1000 $ testing pairs were used. \Cref{fig:gridsweep_darcy1} demonstrates resolution-invariant test error for various $ m $, while \Cref{fig:gridsweep_darcy2} displays the relative error of the learned coefficient $ \al^{(r)} $ at resolution $ r $ w.r.t. the coefficient learned on the highest resolution ($ r=129 $).}
	\label{fig:gridsweep_darcy}
\end{figure}

In~\Cref{fig:gridsweep_darcy}, we again confirm that our method is invariant to the refinement of the mesh and improves with more random features. While the difference at low resolutions is more pronounced than that observed for Burgers' equation, our results for Darcy flow still suggest that the expected relative test error converges to a constant value as resolution increases; an estimate of this rate of convergence is seen in~\Cref{fig:gridsweep_darcy2}, where we plot the relative error of the learned parameter $ \al^{(r)} $ at resolution $ r $ w.r.t. the parameter learned at the highest resolution trained, which was $ r=129 $. Although we do not observe the limiting error following the Monte Carlo rate in $ m $, which suggests that the RKHS $ \cH_{k_{\mu}} $ induced by the choice of $ \varphi $ may not be expressive enough (e.g., not universal \cite{sun2018random}), the numerical results make clear that our method nonetheless performs well as an operator approximator.

\section{Conclusions}
\label{sec:conclusion}
In this article, we introduced a random feature methodology for the data-driven approximation of maps between infinite-dimensional Banach spaces. The random feature model, as an emulator of such maps, performs dimension reduction in the sense that the original infinite-dimensional learning problem reduces to an approximate problem of finding $ m $ real numbers (\Cref{sec:problem}). 
Our conceptually infinite-dimensional algorithm is non-intrusive and
results in a scalable method that is consistent with the continuum limit, robust to discretization, and highly flexible in practical use. These benefits were verified in numerical experiments for two nonlinear forward operators based on PDEs, one involving a semigroup and another a coefficient-to-solution operator~(\Cref{sec:experiment}).
While the random feature-based operator emulator learned from data is not guaranteed to be cheaper to evaluate than a full order solver in general, our design of problem-specific random feature maps in~\Cref{sec:application} leads to efficient $ O(mK\log K) $ evaluation of an $ m $-term RFM for simple physical domain geometries and hence competitive computational cost in many-query settings. A straightforward GPU implementation would provide further acceleration.

There are various directions for future work. We are interested in application of random feature methods to more challenging problems in the sciences, such as climate modeling and material modeling, and to the solution of design and inverse problems arising in those settings with the RFM serving as a cheap emulator. Of great importance in furthering the methodology is the question of how to adapt the random features to data instead of manually constructing them. Some possibilities along this line of work include the Bayesian optimization of RFM hyperparameters, as frequently used in Gaussian process regression, or more general hierarchical learning of the pair $ (\varphi, \mu) $ itself, both of which would lead to data-adapted induced kernels. Such developments would make the RFM more streamlined, competitive with deep learning alternatives, and serve to further clarify the effectiveness of function space learning algorithms. Finally, the development of a theory which underpins our method, allows for proof of convergence, and characterizes the quality of the RKHS spaces induced by random feature maps, would be both mathematically interesting and highly desirable as it 
would help guide methodological development.

\appendix
\section{Proofs of Results}\label{sec:appendix_proofs}
\begin{proof}[Proof of \cref{r:1}]\label{proof:result1}
	Fix $ a\in\cX $ and $ y\in\cY $. Then, we note that
	\begin{equation}\label{eqn:rkhs_proof_prop1}
	k_{\mu}(\cdot,a)y=\int\ip{\varphi(a;\theta)}{y}_{\cY}\varphi(\cdot;\theta)\mu(d\theta)=\cA\ip{\varphi(a;\cdot)}{y}_{\cY}\in\im(\cA)\, ,
	\end{equation}
	since $ \ip{\varphi(a;\cdot)}{y}_{\cY}\in L_{\mu}^{2}(\Theta;\R) $ by the Cauchy-Schwarz inequality.
	
	Now we show that $ \im(\cA) $ admits a reproducing property of the form~\cref{eqn:rkprop_banach}. First, note that $ \cA $ can be viewed as a bijection between its coimage and image spaces, and we denote this
	bijection by
	\begin{equation}\label{eqn:bijection}
	\tilde{\cA}: \ker(\cA)^{\perp}\to \im(\cA)\, .
	\end{equation}
	For any $ F,\, G \in\im(\cA)$, define the candidate RKHS inner product $ \ip{\cdot}{\cdot} $ by
	\begin{equation}\label{eqn:ip_candidate}
	\ip{F}{G}\defby \ip[\big]{\tilde{\cA}^{-1}F}{\tilde{\cA}^{-1}G}_{L^{2}_{\mu}(\Theta;\R)}\, .
	\end{equation}
	This is indeed a valid inner product since $ \tilde{\cA} $ is invertible. Note that for any $ q\in\ker(\cA) $,
	\begin{align*}
	\ip[\big]{q}{\ip{\varphi(a;\cdot)}{y}_{\cY}}_{L^{2}_{\mu}(\Theta;\R)}&=\int q(\theta)\ip{\varphi(a;\theta)}{y}_{\cY}\, \mu(d\theta)\\
	&=\ip[\Big]{\int q(\theta)\varphi(a;\theta)\mu(d\theta)}{y}_{\cY}\\
	&=0
	\end{align*}
	so that $ \ip{\varphi(a;\cdot)}{y}_{\cY}\in \ker(\cA)^{\perp} $. Then for any $ F\in \im(\cA) $, we compute
	\begin{align*}
	\ip{k_{\mu}(\cdot,a)y}{F}&=\ip[\big]{\ip{\varphi(a;\cdot)}{y}_{\cY}}{\tilde{\cA}^{-1}F}_{L^{2}_{\mu}(\Theta;\R)}\\
	&=\int \ip{\varphi(a;\theta)}{y}_{\cY}(\tilde{\cA}^{-1}F)(\theta)\mu(d\theta)\\
	&=\ip[\Big]{\int(\tilde{\cA}^{-1}F)(\theta)\varphi(a;\theta)\mu(d\theta)}{y}_{\cY}\\
	&=\ip[\big]{y}{(\cA\tilde{\cA}^{-1}F)(a)}_{\cY}\\
	&=\ip{y}{F(a)}_{\cY}\, ,
	\end{align*}
	which gives exactly~\cref{eqn:rkprop_banach} if our
	candidate inner product is defined to be the RKHS inner product. 
	Since $ F\in\im(\cA) $ is arbitrary, this and~\cref{eqn:rkhs_proof_prop1} 
	together imply that $ \im(\cA)=\cH_{k_{\mu}} $ is the RKHS induced by $ k_{\mu} $ as shown in~\cite{cucker2002mathematical,kadri2016operator}. 
\end{proof}

\begin{proof}[Proof of \cref{r:2}]
	Since $ L_{\mu^{(m)}}^{2}(\Theta;\R) $ is isomorphic to $ \R^{m} $, we can consider the map $ \cA: \R^{m}\to L^{2}_{\nu}(\cX;\cY) $ defined in~\cref{eqn:rf_integral_operator} and use \cref{r:1} to conclude that
	\begin{equation}\label{eqn:rkhs_finite_dim}
	\cH_{k^{(m)}}=\im(\cA)=\biggl\{ \dfrac{1}{m}\sum_{j=1}^{m}c_j\varphi(\cdot;\theta_{j}): c\in\R^{m} \biggr\} = \spanvec\{\varphi_{j}\}_{j=1}^{m}\, ,
	\end{equation}
	since the $ \{\varphi_{j}\}_{j=1}^{m} $ are assumed linearly independent.
\end{proof}

\begin{proof}[Proof of \cref{res:minimizer_equiv}]
	Recall from \cref{r:2} that the RKHS $\cH_{k^{(m)}}$ comprises the linear span of the $\{\varphi_j\defby \varphi(\cdot;\theta_j)\}_{j=1}^m$. 
	Hence $ \varphi_{j}\in\cH_{k^{(m)}} $, and note that by the reproducing kernel property \cref{eqn:rkprop_banach}, for any $F \in \cH_{k^{(m)}}$, $a\in \cX$ and $y \in \cY$,
	\begin{align*}
	\ip{y}{F(a)}_{\cY}&=\ip[\big]{k^{(m)}(\cdot,a)y}{F}_{\cH_{k^{(m)}}}\\
	&=\frac{1}{m}\sum_{j=1}^m\ip{\varphi_j(a)}{y}_{\cY}\ip{\varphi_j}{F}_{\cH_{k^{(m)}}}\\
	&=\ip[\bigg]{y}{\frac{1}{m}\sum_{j=1}^m \ip{\varphi_j}{F}_{\cH_{k^{(m)}}} \varphi_j(a)}_{\cY}\, .
	\end{align*}
	Since this is true for all $y \in \cY$, we deduce that
	\begin{equation}\label{eqn:rkhs_fd_unique_rep}
	F=\frac{1}{m}\sum_{j=1}^m \alpha_j \varphi_j\, ,\quad\alpha_j=\ip{\varphi_j}{F}_{\cH_{k^{(m)}}}\, .
	\end{equation}
	As the $\{\varphi_j\}_{j=1}^m$ are assumed linearly independent, we deduce that the representation~\cref{eqn:rkhs_fd_unique_rep} is unique.
	
	Finally, we calculate the RKHS norm of any such $F$ in terms of $\alpha$:
	\begin{align*}
	\norm*{F}_{\cH_{k^{(m)}}}^2=\ip*{F}{F}_{\cH_{k^{(m)}}}
	&=\ip[\bigg]{\frac{1}{m}\sum_{j=1}^m \alpha_j \varphi_j}{F}_{\cH_{k^{(m)}}}\\
	&=\frac{1}{m}\sum_{j=1}^m \alpha_j\ip{\varphi_j}{F}_{\cH_{k^{(m)}}}\\
	&=\frac{1}{m}\sum_{j=1}^m \alpha_j^2\, .
	\end{align*}
	Substituting this into~\cref{eqn:opt_equivalence}, we obtain the desired equivalence with \cref{eqn:opt_randomfeature}. 
\end{proof}

\section{Further Remarks on Integral Representation of RKHS}\label{sec:integral_rkhs}
	We recall the linear operator $ \cA $~\cref{eqn:rf_integral_operator} from \Cref{sec:rfm}. In this appendix, we clarify the meaning of \Cref{eqn:rfm_coeff_func} and show that $ \cA $ is a square root of $ T_{k_{\mu}} $. Similar discussion is provided by Bach in \cite{bach2017equivalence}, Sec. 2, for the special case $ \cY=\R $. 
	
	By the assumption $ \varphi\in L_{\nu\times\mu}^2(\cX\times\Theta;\cY) $ and Cauchy-Schwarz inequality, we find
	\begin{equation}\label{eqn:bounded_A}
	\cA\in\cL\bigl(L^2_{\mu}(\Theta;\R), L_{\nu}^2(\cX;\cY)\bigr) \,.
	\end{equation}
	Now let $ F \in \im(\cA)=\cH_{k_{\mu}} $. We have $ F=\cA c $ for some $ c\in L^2_{\mu} $. But since $ \ker(\cA) $ is closed, $ L^2_{\mu}=\ker(\cA)\oplus \ker(\cA)^{\perp}  $ and hence there exist unique $ q_{F}\in\ker(\cA) $ and $ c_{F}\in\ker(\cA)^{\perp} $ such that $ c=q_{F}+c_{F} $. Using the notation in \cref{eqn:bijection}, we have $ c_F=\tilde{\cA}^{-1}F $ by definition of $ \tilde{\cA} $. The reproducing property in \cref{proof:result1} produced the representation $ F=\cA c_F $; in fact, the similar calculation leading to \cref{eqn:rfm_coeff_func} in \Cref{sec:rfm} also identified the unique $ c_F $, there defined formally by $ c_{F}(\theta) = \ip{\varphi(\cdot;\theta)}{F}_{\cH_{k_{\mu}}} $. Indeed,
	\begin{align*}
	\ip{c_{F}}{q}_{L^{2}_{\mu}(\Theta;\R)}&=\int \ip{\varphi(\cdot;\theta)}{F}_{\cH_{k_{\mu}}} q(\theta)\mu(d\theta)\\
	&=\ip[\Big]{\int q(\theta)\varphi(\cdot;\theta)\mu(d\theta)}{F}_{\cH_{k_{\mu}}}\\
	&=0
	\end{align*}
	for any $ q\in\ker(\cA) $. Hence $ c_{F}\in\ker(\cA)^{\perp} $, and we interpret \Cref{eqn:rfm_coeff_func} as formal notation for the unique element $ \tilde{\cA}^{-1}F\in\ker(\cA)^{\perp} $. Using formula \cref{eqn:ip_candidate} and orthogonality, we also obtain the following useful characterization of the RKHS norm:
	\begin{equation}\label{eqn:rkhs_norm_argmin}
	\norm*{F}_{\cH_{k_{\mu}}}^{2}=\norm[\big]{\tilde{\cA}^{-1}F}^{2}_{L^{2}_{\mu}}=\norm*{c_{F}}^{2}_{L^{2}_{\mu}} = \min_{c\in \mathcal{C}_F}\norm{c}^{2}_{L^2_{\mu}}\, ,
	\end{equation}
	where $ \mathcal{C}_F\defby\{c\in L^2_{\mu}(\Theta;\R): \cA c = F\} $.
	
	Finally, we show that $ \cA\cA^{*}=T_{k_{\mu}} $. This means that the RKHS is equal to the image of two different square roots of integral operator $ T_{k_{\mu}} $: $ \cH_{k_{\mu}}=\im(T_{k_{\mu}}^{1/2})=\im(\cA) $.
	First, for any $ F\in L^{2}_{\nu}(\cX;\cY)$ and $ c\in L_{\mu}^{2}(\Theta;\R) $,
	\begin{align*}
	\ip{F}{\cA c}_{L_{\nu}^{2}}&=\ip[\Big]{F}{\int c(\theta)\varphi(\cdot;\theta)\mu(d\theta)}_{L_{\nu}^{2}}\\
	&=\int c(\theta)\ip{F}{\varphi(\cdot;\theta)}_{L_{\nu}^{2}}\,\mu(d\theta)\\
	&=\ip[\Big]{\int \ip{F(a')}{\varphi(a';\cdot)}_{\cY}\,\nu(da')}{c}_{L_{\mu}^{2}}
	\end{align*}
	by the Fubini-Tonelli theorem. So, we deduce that the adjoint of $ \cA $ is
	\begin{align}\label{eqn:rf_integral_operator_adjoint}
	\begin{split}
	\cA^{*}: \ L^2_{\nu}(\cX;\cY) &\to L^{2}_{\mu}(\Theta;\R) \\
	F &\mapsto \cA^{*} F\defby \int \ip{F(a')}{\varphi(a';\cdot)}_{\cY}\,\nu(da') \, ,
	\end{split}
	\end{align}
	which is bounded since $ \cA $ is bounded. For any $ F \in L_{\nu}^2(\cX;\cY) $, we compute
	\begin{align*}
	\cA\cA^{*}F&=\int_{\Theta} (\cA^{*}F)(\theta)\varphi(\cdot;\theta)\mu(d\theta)\\
	&=\int_{\Theta} \int_{\cX} \ip{F(a')}{\varphi(a';\theta)}_{\cY}\, \varphi(\cdot;\theta)\nu(da')\mu(d\theta)\\
	&=\int_{\cX} \left( \int_{\Theta} \varphi(\cdot;\theta)\otimes\varphi(a';\theta)\mu(d\theta) \right) F(a')\nu(da')\\
	&=T_{k_{\mu}}F\, ,
	\end{align*}
	again by Fubini-Tonelli, as desired.

\section*{Acknowledgments}
The authors thank Bamdad Hosseini and Nikola B. Kovachki for helpful discussions and are grateful to the two anonymous referees for their careful reading and insightful comments.

\bibliographystyle{siamplain}
\bibliography{references}

\end{document}